\documentclass[11pt]{article}

\usepackage{amsmath}
\usepackage{amsthm}
\usepackage{amssymb}
\usepackage{amscd}
\usepackage{a4wide}
\usepackage{epic}
\usepackage{eepic}

\theoremstyle{definition}
\newtheorem{theorem}{Theorem}[section]
\newtheorem{prop}[theorem]{Proposition}
\newtheorem{lemma}[theorem]{Lemma}
\newtheorem{corollary}[theorem]{Corollary}
\newtheorem{definition}[theorem]{Definition}
\newtheorem{example}[theorem]{Example}
\newtheorem{remark}[theorem]{Remark}
\newtheorem{conjecture}[theorem]{Conjecture}

\numberwithin{equation}{section}

\newenvironment{demo}[1]{%
  \trivlist
  \item[\hskip\labelsep
        {\it #1.}]
}{%
\hfill\qedsymbol
  \endtrivlist
}

\newcommand\Nat{\mathbb{N}}
\newcommand\Int{\mathbb{Z}}
\newcommand\Rat{\mathbb{Q}}
\newcommand\Real{\mathbb{R}}
\newcommand\Comp{\mathbb{C}}

\newcommand\Par{\operatorname{Par}}
\newcommand\SPar{\operatorname{SPar}}
\newcommand\FF{\mathcal{F}}
\newcommand\GG{\mathcal{G}}
\newcommand\PP{\mathbb{P}}
\newcommand\LL{\mathcal{L}}
\newcommand\Pfun{\mathcal{P}}
\newcommand\Qfun{\mathcal{Q}}
\newcommand\Rfun{\mathcal{R}}

\newcommand\Pf{\operatorname{Pf}}

\newcommand\Symp{\mathbf{Sp}}
\newcommand\twedge{{\textstyle\bigwedge}} 

\newcommand\vectx{\boldsymbol{x}}
\newcommand\vecty{\boldsymbol{y}}
\newcommand\vectz{\boldsymbol{z}}
\newcommand\vecta{\boldsymbol{a}}

\newcommand\trans{{}^t\!}

\newcommand\ep{\varepsilon}
\renewcommand\tilde{\widetilde}

\newcommand\QTab{\operatorname{QTab}}
\newcommand\PTab{\operatorname{PTab}}
\newcommand\wt{\operatorname{wt}}

\allowdisplaybreaks

\title{
Symplectic $Q$-Functions
}

\author{
Soichi Okada%
\footnote{
Graduate School of Mathematics, Nagoya University, 
Furo-cho, Chikusa-ku, Nagoya 464-8602, Japan, 
{\tt okada@math.nagoya-u.ac.jp}
}
\footnote{
This work was partially supported by 
JSPS Grants-in-Aid for Scientific Research No.~24340003, No.~15K13425 and No.~18K03208.
}
}

\date{}

\begin{document}

\maketitle

\begin{abstract}
Symplectic $Q$-functions are a symplectic analogue of Schur $Q$-functions 
and defined as the $t=-1$ specialization of Hall--Littlewood functions 
associated with the root system of type $C$. 
In this paper we prove that symplectic $Q$-functions share many of the properties 
of Schur $Q$-functions, such as a tableau description and a Pieri-type rule.
And we present some positivity conjectures, including the positivity conjecture of structure constants 
for symplectic $P$-functions.
We conclude by giving a tableau description of factorial symplectic $Q$-functions.
\par
Mathematics Subject Classification (MSC2010): 05E05 (primary), 05A15, 05E15 (secondary)
\end{abstract}

\section{%
Introduction
}

The aim of this paper is to establish (and conjecture) properties of symplectic $P$-/$Q$-functions, 
which are a symplectic analogue of Schur $P$-/$Q$-functions 
defined in terms of Hall--Littlewood functions. 
It turns out that symplectic $P$-/$Q$-functions share many of the nice properties of 
Schur $P$-/$Q$-functions.

Schur $Q$-functions, introduced by Schur \cite{Schur}, are an important family of symmetric functions 
as well as Schur ($S$-)functions are.
Schur functions and Schur $Q$-functions appear in similar situations with various applications, 
and enjoy similar properties.
For example, Schur functions describe the characters of irreducible linear representations of symmetric groups, 
and Schur $Q$-functions play the same role for projective representations.
From the combinatorial point of view, both of them are expressed as multivariate generating functions 
of certain tableaux.
And one of the remarkable feature is that they have positive structure constants.

Schur functions and Schur $Q$-functions are obtained from Hall--Littlewood functions 
by specializing the parameter $t$ to $0$ and $-1$ respectively.
Macdonald \cite[\S10]{Macdonald00} introduced a generalization of Hall--Littlewood functions 
to any root system.
In this paper, we study the $t=-1$ specialization of Macdonald's Hall--Littlewood functions 
associated with the root system of type $C$, which we call symplectic $P$-/$Q$-functions.

Let us explain our results in more detail.
A \emph{partition of length $l$} is a weakly decreasing sequence 
$\lambda = (\lambda_1, \dots, \lambda_l)$ of \emph{positive} integers.
We write $l = l(\lambda)$ and $|\lambda| = \sum_{i=1}^l \lambda_i$.

Let $\Phi$ be the root system of type $C_n$ with positive system given by
$$
\Phi^+
 =
\{ \ep_i \pm \ep_j : 1 \le i < j \le n \} \cup \{ 2 \ep_i : 1 \le i \le n \},
$$
where $(\ep_1, \dots, \ep_n)$ is the standard orthonormal basis of $\Real^n$.
We denote by $W_n$ the Weyl group of $\Phi$, which is the semi-direct product $S_n \ltimes (\Int/2\Int)^n$ 
of the symmetric group $S_n$ and an elementary abelian $2$-group $(\Int/2\Int)^n$.
Then $W_n$ acts on the Laurent polynomial ring in $\vectx = (x_1, \dots, x_n)$ 
by permuting and inverting the variables.
We define the \emph{symplectic Hall--Littlewood function} $\PP^C_\lambda(\vectx:t)$ 
corresponding to a partition $\lambda$ of length $l \le n$ by putting
\begin{equation}
\label{eq:HLC}
\PP^C_\lambda (\vectx;t)
 =
\frac{1}{v^{(n)}_\lambda(t)}
\sum_{w \in W_n}
 w \left(
  \prod_{i=1}^n x_i^{\lambda_i}
  \prod_{i=1}^n \frac{ 1 - t x_i^{-2} }{ 1 - x_i^{-2} }
  \prod_{1 \le i < j \le n}
   \frac{ 1 - t x_i^{-1} x_j }{ 1 - x_i^{-1} x_j }
   \frac{ 1 - t x_i^{-1} x_j^{-1} }{ 1 - x_i^{-1} x_j^{-1} }
 \right),
\end{equation}
where $\lambda$ is identified with a sequence $(\lambda_1, \dots, \lambda_n)
 = (\lambda_1, \dots, \lambda_l, 0, \dots, 0)$ of length $n$, 
and the normalizing constant $v^{(n)}_\lambda(t)$ is given explicitly by
\begin{equation}
\label{eq:v}
v^{(n)}_\lambda(t)
 =
\prod_{j=1}^{m_0} \frac{ 1 - t^{2j} }{ 1 - t }
\cdot
\prod_{k \ge 1} \prod_{j=1}^{m_k} \frac{ 1 - t^j }{ 1 - t }
\end{equation}
with $m_k = \# \{ i : 1 \le i \le n, \ m_i = k \}$ ($k \ge 0$).
It can be shown that $\PP^C_\lambda(\vectx;t)$ is a $W_n$-invariant Laurent polynomial 
with coefficients in $\Int[t]$.

The \emph{symplectic Schur function}, denoted by $S^C_\lambda(\vectx)$, 
corresponding to a partition $\lambda$ of length $l \le n$ is defined 
by specializing $t=0$ in the symplectic Hall--Littlewood function:
\begin{equation}
\label{eq:SC=HLC}
S^C_\lambda(\vectx) = \PP^C_\lambda(\vectx;0).
\end{equation}
Then it follows from Weyl's character formula that $S^C_\lambda(\vectx)$ 
is the character of the irreducible highest weight representation 
of the symplectic group $\Symp_{2n}(\Comp)$ with highest weight 
$\lambda_1 \ep_1 + \dots + \lambda_n \ep_n$.
King \cite{King76} introduced the notion of symplectic tableaux 
with entries $1, \overline{1}, \dots, n, \overline{n}$, 
whose weighted generating functions are the symplectic Schur functions.
Since the product of two symplectic Schur functions corresponds to 
the tensor product of associated representations, 
we see that the product is positively expanded as a linear combination 
of symplectic Schur functions.

Now we define symplectic $P$-/$Q$-functions.
A partition $\lambda$ of length $l$ is called \emph{strict} if $\lambda_1 > \dots > \lambda_l$.
Given a strict partition of length $l \le n$, 
we define the corresponding \emph{symplectic $P$-function} $P^C_\lambda(\vectx)$ 
and \emph{symplectic $Q$-function} $Q^C_\lambda(\vectx)$ by putting
\begin{equation}
\label{eq:PC=HLC}
P^C_\lambda(\vectx)
 =
\PP^C_\lambda(\vectx;-1),
\quad
Q^C_\lambda(\vectx)
 =
2^{l(\lambda)} \PP^C_\lambda(\vectx;-1)
\end{equation}
respectively.
Then we can use the theory of generalized $P$-/$Q$-functions developed in \cite{Okada20} 
to derive several formulas for symplectic $P$-/$Q$-functions.

One of our main results is a proof of King--Hamel's conjecture \cite{KH07} 
on a tableau description of symplectic $Q$-functions.
The \emph{shifted diagram} $S(\lambda)$ of a strict partition $\lambda$ is defined by
$$
S(\lambda)
 =
\{ (i,j) \in \Int^2 : 1 \le i \le l(\lambda), \, i \le j \le \lambda_i +i-1 \}
$$
and represented by replacing lattice points with unit cells.
A shifted tableau of shape $\lambda$ is a filling of the cells of $S(\lambda)$.

\begin{theorem}
\label{thm:main1}
(the non-skew case of Theorem~\ref{thm:tableaux}, conjectured in \cite[Conjecture~3.1]{KH07})
For a sequence $\vectx = (x_1, \dots, x_n)$ of indeterminates 
and a strict partitions $\lambda$ of length $\le n$, we have
$$
Q^C_\lambda(\vectx)
 =
\sum_{T \in \QTab^C_n(\lambda)} \vectx^T,
$$
where the sum is taken over all symplectic primed shifted tableaux of shape $\lambda$ 
with entries $1', 1, \overline{1}', \overline{1}, \dots, n', n, \overline{n}', \overline{n}$, 
and entries $i'$ and $i$ (resp. $\overline{i}'$ and $\overline{i}$) in $T$ 
contribute with $x_i$ (resp. $x_i^{-1}$).
See Definition~\ref{def:tableaux} for a precise definition.
\end{theorem}

More generally, we obtain the skew version (Theorem~\ref{thm:tableaux}) 
and the factorial extension (Theorem~\ref{thm:fac_skew_tableaux}) 
of this theorem.

For three strict partitions $\lambda$, $\mu$ and $\nu$ of length $\le n$, 
let $\tilde{f}^\lambda_{\mu,\nu}$ be the structure constant determined by the formula
\begin{equation}
\label{eq:str_const}
P^C_\mu (\vectx) \cdot P^C_\nu (\vectx)
 =
\sum_\lambda \tilde{f}^\lambda_{\mu,\nu} P^C_\lambda(\vectx),
\end{equation}
Another main result of this paper is the Pieri-type rule for symplectic $P$-functions, 
which gives an explicit formula for the structure constants $\tilde{f}^\lambda_{\mu,(r)}$ 
in the case where $\nu = (r)$ is a one-row partition.

\begin{theorem}
\label{thm:main2}
(Theorem~\ref{thm:Pieri})
For two strict partitions $\lambda$, $\mu$ and a positive integer $r$, we have
$$
\tilde{f}^\lambda_{\mu,(r)}
 =
\begin{cases}
 \sum_\kappa 2^{a(\mu,\kappa) + a(\lambda,\kappa) - \chi[l(\mu)>l(\kappa)] - 1} 
 &\text{if $l(\lambda) = l(\mu)$ or $l(\mu)+1$,}
\\
 0 &\text{otherwise,}
\end{cases}
$$
where $\kappa$ runs over all strict partitions satisfying 
$\mu_1 \ge \kappa_1 \ge \mu_2 \ge \kappa_2 \ge \dots$, 
$\lambda_1 \ge \kappa_1 \ge \lambda_2 \ge \kappa_2 \ge \dots$, 
and $(|\mu|-|\kappa|)+(|\lambda-|\kappa|) = r$.
And $a(\mu,\kappa)$ (resp. $a(\lambda,\kappa)$) is the number of connected components of 
the skew shifted diagram $S(\mu) \setminus S(\kappa)$ (resp. $S(\lambda) \setminus S(\kappa)$), 
and $\chi[l(\mu)>l(\kappa)] = 1$ if $l(\mu) > l(\kappa)$ and $0$ otherwise.
\end{theorem}

Based on this result and computer experiments, we propose the following conjecture.

\begin{conjecture}
\label{conj:main3}
(Conjecture~\ref{conj:1})
The structure constants $\tilde{f}^\lambda_{\mu,\nu}$ defined by (\ref{eq:str_const}) are 
nonnegative integers.
\end{conjecture}

Moreover we give several positivity conjectures on various expansion coefficients 
involving symplectic $P$-/$Q$-functions 
(Conjectures~\ref{conj:2}, \ref{conj:3-1}, \ref{conj:3-2} and \ref{conj:4}).

Here we mention orthogonal $P$-/$Q$-functions, which are defined as the $t=-1$ specialization 
of Hall--Littlewood functions associated with the root system of types $B$ and $D$.
By \cite[Theorem~7.2]{Okada20}, we can express them as generalized $P$-functions associated 
to certain polynomial sequences, so we obtain the Nimmo-type formula like (\ref{eq:Nimmo_QC}) 
and the Schur-type Pfaffian formula like (\ref{eq:Schur_QC}).
However orthogonal $P$-/$Q$-functions do not behave as nicely as symplectic $P$-/$Q$-functions.
For example, some formulas take different forms depending the parity of the rank of the root system, 
and some structure constants with respect to the odd orthogonal $P$-functions are negative. 

The remainder of this paper is organized as follows.
In Section~2, we apply results on generalized $P$-functions given in \cite{Okada20} 
to obtain several properties of symplectic $P$-/$Q$-functions.
In Section~3, we lift symplectic $P$-/$Q$-functions from Laurent polynomials in finitely many variables 
to symmetric functions in infinitely many variables, and introduce the notion of skew symplectic $Q$-functions.
Sections~4 and 5 are devoted to the proof of the tableau description and the Pieri-type rule respectively.
In Section~6, we present some positivity conjectures on various expansion coefficients 
involving symplectic $P$-functions. 
In Section~7, we introduce a factorial analogue of symplectic $Q$-functions 
as the $t=-1$ specialization of Hall--Littlewood-type function 
and give a tableau description for them.
\section{%
Basic properties of symplectic $P$-/$Q$-functions
}

In this section, we collect several properties of symplectic $P$-/$Q$-functions, 
which will be used in our discussion below.
These properties follow from results on generalized $P$-functions proved in \cite{Okada20}.

Recall the definition of generalized $P$-functions associated with 
a polynomial sequence (\cite[Definition~1.2]{Okada20}).
Let $\FF = \{ f_d(u) \}_{d=0}^\infty$ be a sequence of polynomials in one variable $u$ 
such that $f_0(u) = 1$ and $\det f_d(u) = d$ for $d \ge 0$.
Given a sequence $\vectx = (x_1, \dots, x_n)$ of $n$ indeterminates 
and a strict partition $\lambda$ of length $l$, 
we define the \emph{generalized $P$-function} $P^{\FF}_\lambda(\vectx)$ 
associated with $\FF$ by
\begin{equation}
\label{eq:genP}
P^\FF_\lambda(\vectx)
 =
\begin{cases}
 \dfrac{ 1 }{ \Delta(\vectx) }
 \Pf \begin{pmatrix}
  A(\vectx) & V^\FF_\lambda(\vectx) \\
  - \trans V^\FF_\lambda(\vectx) & O
 \end{pmatrix}
 &\text{if $n+l$ is even,}
\\
 \dfrac{ 1 }{ \Delta(\vectx) }
 \Pf \begin{pmatrix}
  A(\vectx) & V^\FF_{\lambda^0}(\vectx) \\
  - \trans V^\FF_{\lambda^0}(\vectx) & O
 \end{pmatrix}
 &\text{if $n+l$ is odd,} \\
\end{cases}
\end{equation}
where $\lambda = (\lambda_1, \dots, \lambda_l)$, $\lambda^0 = (\lambda_1, \dots, \lambda_l, 0)$, and
\begin{equation}
\label{eq:A}
A(\vectx) = \left( \frac{ x_j - x_i }{ x_j + x_i } \right)_{1 \le i, j \le n},
\quad
\Delta(\vectx) = \prod_{1 \le i < j \le n} \frac{ x_j - x_i }{ x_j + x_i },
\quad
V^{\FF}_\alpha(\vectx)
 = 
\Big( f_{\alpha_j}(x_i) \Big)_{1 \le i \le n, \, 1 \le j \le r}.
\end{equation}
By the Nimmo formula \cite[(A13)]{Nimmo90}, 
we see that, if the polynomial sequence $\FF$ is given by $f_d(u) = u^d$ (resp. $2 u^d$) for $d \ge 1$, 
then the associated generalized $P$-functions are the classical Schur $P$-functions (resp. $Q$-functions).

The symplectic $P$-/$Q$-functions given by (\ref{eq:PC=HLC}) are written as generalized $P$-functions 
up to a simple transformation of variables.
For a nonnegative integer $d$, we define Laurent polynomials $\tilde{f}_d(x)$, $\tilde{g}_d(x)$ 
and polynomial $f_d(u)$, $g_d(u)$ by
\begin{align}
\label{eq:f}
\tilde{f}_d(x)
 &=
f_d(x+x^{-1})
 =
\begin{cases}
 1 &\text{if $d = 0$,} \\
 (x^d - x^{-d}) \dfrac{ x + x^{-1} }{ x - x^{-1} } &\text{if $d \ge 1$,}
\end{cases}
\\
\label{eq:g}
\tilde{g}_d(x)
 &=
g_d(x+x^{-1})
 =
\begin{cases}
 1 &\text{if $d = 0$,} \\
 2 (x^d - x^{-d}) \dfrac{ x + x^{-1} }{ x - x^{-1} } &\text{if $d \ge 1$.}
\end{cases}
\end{align}
Then we have 

\begin{prop}
\label{prop:PC=genP}
(\cite[Theorem~7.2]{Okada20})
Let $\FF = \{ f_d(u) \}_{d=0}^\infty$ and $\GG = \{ g_d(u) \}_{d=0}^\infty$ be 
the polynomial sequences defined by (\ref{eq:f}) and (\ref{eq:g}) respectively.
For a sequence $\vectx = (x_1, \dots, x_n)$ of variables and a strict partition $\lambda$ of length $\le n$, 
we have
\begin{equation}
\label{eq:PC=genP}
P^C_\lambda(\vectx)
 =
P^{\FF}_\lambda(\vectx+\vectx^{-1}),
\quad
Q^C_\lambda(\vectx)
 =
P^{\GG}_\lambda(\vectx+\vectx^{-1}),
\end{equation}
where $\vectx + \vectx^{-1} = (x_1+x_1^{-1}, \dots, x_n+x_n^{-1})$.
\end{prop}

In particular, if $\vectx = (x)$ is a single variable, then we have
\begin{equation}
\label{eq:n=1}
P^C_{(r)}(x) = \tilde{f}_r(x),
\quad
Q^C_{(r)}(x) = \tilde{g}_r(x).
\end{equation}
By combining Proposition~\ref{prop:PC=genP} with the definition (\ref{eq:genP}) of generalized $P$-functions, 
we obtain the following Nimmo-type expression of symplectic $P$-/$Q$-functions.

\begin{prop}
\label{prop:Nimmo}
Let $\vectx = (x_1, \dots, x_n)$ and put
\begin{equation}
\label{eq:AC}
\tilde{A}(\vectx)
 =
\left(
 \frac{ (x_j + x_j^{-1}) - (x_i + x_i^{-1}) }
      { (x_j + x_j^{-1}) + (x_i + x_i^{-1}) }
\right)_{1 \le i, j \le n},
\quad
\tilde{\Delta}(\vectx)
 =
\prod_{1 \le i < j \le n}
 \frac{ (x_j + x_j^{-1}) - (x_i + x_i^{-1}) }
      { (x_j + x_j^{-1}) + (x_i + x_i^{-1}) },
\end{equation}
and
$$
\tilde{V}_{(\alpha_1, \dots, \alpha_r)}(\vectx)
 = 
\Big( \tilde{f}_{\alpha_j}(x_i) \Big)_{1 \le i \le n, \, 1 \le j \le r},
\quad
\tilde{W}_{(\alpha_1, \dots, \alpha_r)}(\vectx)
 = 
\Big( \tilde{g}_{\alpha_j}(x_i) \Big)_{1 \le i \le n, \, 1 \le j \le r}.
$$
Then, for a strict partition $\lambda$ of length $l \le n$, we have
\begin{align}
\label{eq:Nimmo_PC}
P^C_\lambda(\vectx)
 &=
\begin{cases}
 \dfrac{ 1 }{ \tilde{\Delta}(\vectx) }
 \Pf \begin{pmatrix}
  \tilde{A}(\vectx) & \tilde{V}_\lambda(\vectx) \\
  - \trans \tilde{V}_\lambda(\vectx) & O
 \end{pmatrix}
 &\text{if $n+l$ is even,}
\\
 \dfrac{ 1 }{ \tilde{\Delta}(\vectx) }
 \Pf \begin{pmatrix}
  \tilde{A}(\vectx) & \tilde{V}_{\lambda^0}(\vectx) \\
  - \trans \tilde{V}_{\lambda^0}(\vectx) & O
 \end{pmatrix}
 &\text{if $n+l$ is odd,} \\
\end{cases}
\\
\label{eq:Nimmo_QC}
Q^C_\lambda(\vectx)
 &=
\begin{cases}
 \dfrac{ 1 }{ \tilde{\Delta}(\vectx) }
 \Pf \begin{pmatrix}
  \tilde{A}(\vectx) & \tilde{W}_\lambda(\vectx) \\
  - \trans \tilde{W}_\lambda(\vectx) & O
 \end{pmatrix}
 &\text{if $n+l$ is even,}
\\
 \dfrac{ 1 }{ \tilde{\Delta}(\vectx) }
 \Pf \begin{pmatrix}
  \tilde{A}(\vectx) & \tilde{W}_{\lambda^0}(\vectx) \\
  - \trans \tilde{W}_{\lambda^0}(\vectx) & O
 \end{pmatrix}
 &\text{if $n+l$ is odd,} \\
\end{cases}
\end{align}
where $\lambda^0 = (\lambda_1, \dots, \lambda_l, 0)$.
\end{prop}

In what follows, 
we adopt (\ref{eq:Nimmo_PC}) and (\ref{eq:Nimmo_QC}) 
as the definitions of $P^C_\lambda(\vectx)$ and $Q^C_\lambda(\vectx)$ 
for a general strict partition $\lambda$ respectively.

Recall the following relation between Pfaffians and determinants 
(see e.g. \cite[Corollary~2.4 (1)]{Okada19}): 
\begin{equation}
\label{eq:Pf-Laplace}
\Pf \begin{pmatrix}
 Z & W \\
 -\trans W & O
\end{pmatrix}
 =
\begin{cases}
 (-1)^{n(n-1)/2} \det W &\text{if $n=m$,} \\
 0 &\text{if $n<m$,}
\end{cases}
\end{equation}
where $Z$ is an $n \times n$ skew-symmetric matrix and $W$ is an $n \times m$ matrix.
By using this relation (\ref{eq:Pf-Laplace}), we obtain

\begin{lemma}
\label{lem:PC=0}
Let $\vectx = (x_1, \dots, x_n)$.
Under the definition (\ref{eq:Nimmo_PC}) and (\ref{eq:Nimmo_QC}), 
we have $P^C_\lambda(\vectx) = Q^C_\lambda(\vectx) = 0$ 
for a strict partition $\lambda$ with $l(\lambda) > n$.
\end{lemma}

By \cite[Theorem~2.6]{Okada20} we obtain the following Schur-type Pfaffian expression 
of $Q^C_\lambda(\vectx)$.

\begin{prop}
\label{prop:Schur}
For an arbitrary strict partition $\lambda$, we have
\begin{equation}
\label{eq:Schur_QC}
Q^C_\lambda(\vectx)
 = 
\Pf \Big(
 Q^C_{(\lambda_i,\lambda_j)}(\vectx)
\Big)_{1 \le i, j \le m},
\end{equation}
where $m = l(\lambda)$ or $l(\lambda)+1$ according whether $l(\lambda)$ is even or odd, 
and we use the convention
\begin{equation}
\label{eq:convention}
Q^C_{(s,r)}(\vectx) = - Q^C_{(r,s)}(\vectx),
\quad
Q^C_{(r,0)}(\vectx) = - Q^C_{(0,r)}(\vectx) = Q^C_{(r)}(\vectx),
\quad
Q^C_{(0,0)}(\vectx) = 0
\end{equation}
for positive integers $r$ and $s$.
\end{prop}

The entries of the Pfaffian of (\ref{eq:Schur_QC}) are obtained 
from the following generating functions.

\begin{prop}
(\cite[Proposition~7.6]{Okada20})
\label{prop:length12}
If we put
$$
\tilde{\Pi}_z(\vectx)
 =
\prod_{i=1}^n
 \frac{ (1-x_iz) (1-x_i^{-1}z) }
      { (1+x_iz) (1+x_i^{-1}z) },
$$
then we have
\begin{gather}
\label{eq:GF_length1}
\sum_{r \ge 0} Q^C_{(r)}(\vectx) z^r
 =
\tilde{\Pi}_z(\vectx),
\\
\label{eq:GF_length2}
\sum_{r,s \ge 0} Q^C_{(r,s)}(\vectx) z^r w^s
 =
\frac{ (z-w)(1-zw) }
     { (z+w)(1+zw) }
\left(
 \tilde{\Pi}_z(\vectx) \tilde{\Pi}_w(\vectx) - 1
\right),
\end{gather}
where $Q^C_{(0)}(\vectx) = 1$ and we use the convention (\ref{eq:convention}).
\end{prop}

We can use these generating functions to derive a formula 
for symplectic $Q$-functions for length $2$ partitions 
in terms of those for length $1$ partitions.

\begin{corollary}
\label{cor:length2}
For a strict partition $(r,s)$ of length $2$, we have
\begin{align}
\label{eq:length2}
Q^C_{(r,s)}(\vectx)
&
 =
Q^C_{(r)}(\vectx) Q^C_{(s)} (\vectx)
\notag
\\
&\quad
 +
2 \sum_{k=1}^s (-1)^k 
 \left(
  Q^C_{(r+k)}(\vectx)
  + 2 \sum_{i=1}^{k-1} Q^C_{(r+k-2i)}(\vectx)
  + Q^C_{(r-k)}(\vectx)
 \right)
 Q^C_{(s-k)}(\vectx).
\end{align}
\end{corollary}

\begin{demo}{Proof}
In the proof we simply write $Q_\lambda$ for $Q^C_\lambda(\vectx)$ 
and $Q'_{(r,s)}$ for the right hand side of (\ref{eq:length2}).
We prove $Q_{(r,s)} = Q'_{(r,s)}$ by induction on $s$.
It follows from (\ref{eq:GF_length1}) and (\ref{eq:GF_length2}) that
\begin{equation}
\label{eq:GF_length2a}
(z+w)(1+zw)
\sum_{r,s \ge 0} Q_{(r,s)} z^r w^s
 =
(z-w)(1-zw)
\left(
 \sum_{r \ge 0 } Q_{(r)} z^r \cdot \sum_{s \ge 0} Q_{(s)} w^s - 1
\right).
\end{equation}

Suppose that $r \ge 2$.
By equating the coefficients of $z^{r+1} w$ in (\ref{eq:GF_length2a}), we have
$$
Q_{(r,1)} + Q_{(r+1,0)} + Q_{(r-1,0)}
 =
Q_{(r)} Q_{(1)} - Q_{(r+1)} - Q_{(r-1)},
$$
which implies $Q_{(r,1)} = Q_{(r)} Q_{(1)} - 2 Q_{(r+1)} - 2 Q_{(r-1)} = Q'_{(r,1)}$. 

Suppose that $r > s \ge 2$.
By equating the coefficients of $z^{r+1} w^s$ in (\ref{eq:GF_length2a}), we have
\begin{multline}
\label{eq:recursion1}
Q_{(r,s)} + Q_{(r+1,s-1)} + Q_{(r-1,s-1)} + Q_{(r,s-2)}
\\
 =
Q_{(r)} Q_{(s)} - Q_{(r+1)} Q_{(s-1)} - Q_{(r-1)} Q_{(s-1)} + Q_{(r)} Q_{(s-2)}.
\end{multline}
By the definition of $Q'_{(u,v)}$, we have
\begin{align*}
Q'_{(u,v)} + Q'_{(u-1,v-1)}
 &=
Q_{(u)} Q_{(v)} + 2 \sum_{i=1}^v (-1)^i Q_{(u+i)} Q_{(v-i)}
\\
 &\quad\quad
- Q_{(u-1)} Q_{(v-1)} - 2 \sum_{i=1}^{v-1} (-1)^i Q_{(u-1+i)} Q_{(v-1-i)}.
\end{align*}
Hence, by taking $(u,v) = (r,s)$ and $(r+1,s-1)$, we see that
\begin{multline}
\label{eq:recursion2}
Q'_{(r,s)} + Q'_{(r+1,s-1)} + Q'_{(r-1,s-1)} + Q'_{(r,s-2)}
\\
 =
Q_{(r)} Q_{(s)} - Q_{(r+1)} Q_{(s-1)} - Q_{(r-1)} Q_{(s-1)} + Q_{(r)} Q_{(s-2)}.
\end{multline}
Since $Q_{(r+1,s-1)} = Q'_{(r+1,s-1)}$, $Q_{(r-1,s-1)} = Q'_{(r-1,s-1)}$ and $Q_{(r,s-2)} = Q'_{(r,s-2)}$ 
by the induction hypothesis, 
we can conclude $Q_{(r,s)} = Q'_{(r,s)}$ by comparing (\ref{eq:recursion1}) with (\ref{eq:recursion2}).
\end{demo}

Let $\tilde{\Lambda}^{(n)} = \Rat[x_1^{\pm 1}, \dots, x_n^{\pm 1}]^{W_n}$ be 
the ring of $W_n$-invariant Laurent polynomials in $x_1, \dots, x_n$.
Let $\tilde{\Gamma}^{(n)}$ be the subring of $\tilde{\Lambda}^{(n)}$ 
defined by
$$
\tilde{\Gamma}^{(n)}
 =
\{ 
 g \in \tilde{\Lambda}^{(n)} : 
 \text{$g(t, -t, x_3, \dots, x_n)$ is independent of $t$}
\}.
$$

\begin{prop}
\label{prop:PC_basis}
We denote by $\SPar^{(n)}$ the set of all strict partitions of length $\le n$.
The symplectic $P$-functions $\{ P^C_\lambda(\vectx) : \lambda \in \SPar^{(n)} \}$ 
form a basis of $\tilde{\Gamma}^{(n)}$, 
and so does the symplectic $Q$-functions \{$ Q^C_\lambda(\vectx) : \lambda \in \SPar^{(n)} \}$.
\end{prop}

\begin{demo}{Proof}
Let $\Gamma^{(n)}$ be the subring of $\Lambda^{(n)} = \Rat[z_1, \dots, z_n]^{S_n}$, 
the ring of symmetric polynomials, consisting of $f \in \Lambda^{(n)}$ 
such that $f(u, -u, z_3, \dots, z_n)$ is independent of $u$.
Let $\FF = \{ f_d(u) \}_{d \ge 0}$ be the polynomial sequence given by (\ref{eq:f}).
Then, by \cite[Corollary~5.2]{Okada20}, 
the generalized $P$-functions $\{ P^{\FF}_\lambda(\vectz) : \lambda \in \SPar^{(n)} \}$ 
form a basis of $\Gamma^{(n)}$.

Let $\phi : \Rat[z_1, \dots, z_n] \to \Rat[x_1^{\pm 1}, \dots, x_n^{\pm 1}]$ 
be the ring homomorphism given by $\phi(z_i) = x_i + x_i^{-1}$ ($1 \le i \le n$).
Since $\tilde{\Lambda}^{(n)} = \Rat[x_1+x_1^{-1}, \dots, x_n+x_n^{-1}]^{S_n}$, 
the homomorphism $\phi$ induces an isomorphism $\Lambda^{(n)} \to \tilde{\Lambda}^{(n)}$.
Given a (Laurent) polynomial $h(x_1, \dots, x_n)$, we see that 
$h(t,-t,x_3,\dots, x_n)$ is independent of $t$ if and only if 
$(\partial h/\partial x_1)(t,-t,x_3,\dots,x_n) = 
(\partial h/\partial x_2)(t,-t,x_3,\dots,x_n)$.
By using this observation, we can check that $f \in \Gamma^{(n)}$ if and only if $\phi(f) \in \tilde{\Gamma}^{(n)}$.

Hence the proposition follows from $\phi(P^{\FF}_\lambda(\vectz)) = P^C_\lambda(\vectx)$ 
(Proposition~\ref{prop:PC=genP}).
\end{demo}

We conclude this section with a relation between symplectic $P$-functions and symplectic Schur functions.
For a partition $\mu$ of length $\le n$, 
the definition (\ref{eq:SC=HLC}) of the symplectic Schur function $S^C_\mu (\vectx)$ 
can be rewritten in the bialternant form
\begin{equation}
\label{eq:SC}
S^C_\mu(\vectx)
 =
\frac{ \det \Big( x_i^{\mu_j+n-j+1} - x_i^{-(\mu_j+n-j+1) } \Big)_{1 \le i, j \le n} }
     { \det \Big( x_i^{n-j+1} - x_i^{-(n-j+1) } \Big)_{1 \le i, j \le n} }.
\end{equation}
Note that the denominator factors as
\begin{equation}
\label{eq:SC_denom}
\det \Big( x_i^{n-j+1} - x_i^{-(n-j+1) } \Big)_{1 \le i, j \le n}
 =
\prod_{i=1}^n (x_i - x_i^{-1})
\prod_{1 \le i < j \le n} \big( (x_i+x_i^{-1}) - (x_j+x_j^{-1}) \big).
\end{equation}
Let $\delta_n = (n, n-1, \dots, 2,1)$ be the staircase strict partition of length $n$.

\begin{prop}
\label{prop:PC=SC*SC}
For a partition $\mu$ of length $\le n$, the symplectic $P$-function corresponding to 
a strict partition $\mu+\delta_n = (\mu_1+n, \dots, \mu_n+1)$ can be expressed as the product 
of two symplectic Schur functions:
\begin{equation}
\label{eq:PC=SC*SC}
P^C_{\mu+\delta_n}(\vectx)
 =
S^C_{\delta_n}(\vectx) \cdot S^C_\mu(\vectx).
\end{equation}
\end{prop}

\begin{demo}{Proof}
By applying the formula (\ref{eq:Pf-Laplace}) 
to the Pfaffian in the Nimmo-type formula (\ref{eq:Nimmo_PC}), we have
\begin{align*}
P^C_{\mu+\delta_n}(\vectx)
 &=
\frac{ 1 }
     { \tilde{\Delta}(\vectx) }
\cdot
(-1)^{\binom{n}{2}}
\det \left(
 \big( x_i^{\mu_j+n-j+1} - x_i^{-(\mu_j+n-j+1)} \big) \frac{ x_i + x_i^{-1} }{ x_i - x_i^{-1} }
\right)_{1 \le i, j \le n}
\\
 &=
\frac{ \prod_{i=1}^n (x_i+x_i^{-1}) \prod_{1 \le i < j \le n} \big( (x_i+x_i^{-1}) + (x_j+x_j^{-1}) \big) }
     { \prod_{i=1}^n (x_i-x_i^{-1}) \prod_{1 \le i < j \le n} \big( (x_i-x_i^{-1}) - (x_j+x_j^{-1}) \big) }
\\
&\quad\quad
\times
\det \left(
 x_i^{\mu_j+n-j+1} - x_i^{-(\mu_j+n-j+1)}
\right)_{1 \le i, j \le n}
\end{align*}
Since $S^C_{\delta_n}(\vectx) 
 = \prod_{i=1}^n (x_i+x_i^{-1}) \prod_{1 \le i < j \le n} \big( (x_i+x_i^{-1}) + (x_j+x_j^{-1}) \big)$ 
by (\ref{eq:SC}) and (\ref{eq:SC_denom}), 
we can obtain the desired identity (\ref{eq:PC=SC*SC}).
\end{demo}

\section{%
Universal symplectic $P$-functions
}

In the first half of this section, 
we define a family of symmetric functions, called universal symplectic $P$-/$Q$-functions, in infinitely many variables, 
which are a lift of symplectic $P$-/$Q$-functions in finitely many variables.
And in the second half we introduce a skew version of universal symplectic $Q$-functions.

\subsection{%
Universal symplectic $P$-functions
}

We begin with introducing Schur $Q$-functions as symmetric functions.
We refer the reader to \cite[III.8]{Macdonald95} for an exposition of Schur $Q$-functions. 
Let $\Lambda = \Lambda(X)$ be the ring of symmetric functions in a countably infinite number of variables $X = \{ x_1, x_2, \dots \}$ 
with rational coefficients.
Let $q_r$ ($r \ge 0$) be the symmetric functions defined by
\begin{equation}
\label{eq:GF_q}
\sum_{r \ge 0} q_r(X) z^r = \prod_{i \ge 1} \frac{1 + x_i z}{1 - x_i z}.
\end{equation}
We denote by $\Gamma$ be the subring of $\Lambda$ generated by $q_r$ ($r \ge 1$), 
and by $\Gamma_d$ the subspace of $\Gamma$ consisting of homogeneous elements of degree $d$.

We define Schur $Q$-functions $Q_\lambda$ as symmetric functions inductively on the length of $\lambda$, 
as Schur \cite[Abschnitt~IX]{Schur} did.
For the empty partition $\emptyset$, we define $Q_\emptyset = 1$, 
and for strict partitions of length $1$ and $2$, we define 
\begin{gather}
\label{eq:Q_length1}
Q_{(r)} = q_r
\quad(r>0),
\\
\label{eq:Q_length2}
Q_{(r,s)}
 =
q_r q_s
 +
2 \sum_{k=1}^s (-1)^k q_{r+k} q_{s-k}
\quad(r>s>0).
\end{gather}
Then the Schur $Q$-function $Q_\lambda$ corresponding to an arbitrary strict partition $\lambda$ is defined by
\begin{equation}
\label{eq:Q}
Q_\lambda
 =
\Pf \Big(
 Q_{(\lambda_i,\lambda_j)}
\Big)_{1 \le i, j \le m}
\end{equation}
where $m = l(\lambda)$ or $l(\lambda)+1$ according whether $l(\lambda)$ is even or odd, and 
we use the convention 
$Q_{(s,r)} = - Q_{(r,s)}$, $Q_{(r,0)} = - Q_{(0,r)} = Q_{(r)}$, 
and $Q_{(0,0)} = 0$ for positive integers $r$ and $s$.
And Schur $P$-functions $P_\lambda$ are given by
\begin{equation}
\label{eq:P}
P_\lambda = 2^{-l(\lambda)} Q_\lambda.
\end{equation}
If we set $x_{n+1} = x_{n+2} = \dots = 0$ in $P_\lambda$ and $Q_\lambda$, 
then we obtain the $t=-1$ specialization of Hall--Littlewood polynomials in $(x_1, \dots, x_n)$.

\begin{prop}
\label{prop:Q_basis}
(see \cite[III (8.9)]{Macdonald95})
Let $\SPar_d$ be the set of all strict partitions of $d$.
Then both $\{ P_\lambda : \lambda \in \SPar_d \}$ and $\{ Q_\lambda : \lambda \in \SPar_d \}$ 
are bases of $\Gamma_d$.
\end{prop}

In a similar way, we introduce another family of symmetric functions, 
which we call universal symplectic $P$-/$Q$-functions.
The \emph{universal symplectic $Q$-functions} $\Qfun^C_\lambda$ are defined 
by induction on the length $l(\lambda)$.
We put $\Qfun^C_\emptyset = 1$, and 
\begin{gather}
\label{eq:uQC_length1}
\Qfun^C_{(r)} = q_r \quad(r > 0),
\\
\label{eq:uQC_length2}
\Qfun^C_{(r,s)}
 =
q_r q_s 
 +
2 \sum_{k=1}^s (-1)^k 
 \left(
  q_{r+k} + 2 \sum_{i=1}^{k-1} q_{r+k-2i} + q_{r-k}
 \right)
 q_{s-k}
\quad(r>s>0).
\end{gather}
Then, for any strict partition $\lambda$, we define
\begin{equation}
\label{eq:uQC}
\Qfun^C_\lambda
 =
\Pf \Big(
 \Qfun^C_{(\lambda_i,\lambda_j)}
\Big)_{1 \le i, j \le m}
\end{equation}
where $m = l(\lambda)$ or $l(\lambda)+1$ according whether $l(\lambda)$ is even or odd, and 
\begin{equation}
\label{eq:univ_convention}
\Qfun^C_{(r,s)} = - \Qfun^C_{(s,r)},
\quad 
\Qfun^C_{(r,0)} = - \Qfun_{(0,r)} = q_r,
\quad
\Qfun^C_{(0,0)} = 0
\end{equation}
for positive integers $r$ and $s$.
The \emph{universal symplectic $P$-functions} $\Pfun^C_\lambda$ are given by 
\begin{equation}
\label{eq:uPC}
\Pfun^C_\lambda = 2^{-l(\lambda)} \Qfun^C_\lambda.
\end{equation}
Universal symplectic $P$-/$Q$-functions are a lift of symplectic $P$-/$Q$-functions 
to the symmetric functions in the following sense.

\begin{prop}
\label{prop:universal}
Let $\tilde{\pi}_n : \Lambda \to \tilde{\Lambda}^{(n)} = \Rat[x_1^{\pm 1}, \dots, x_n^{\pm 1}]^{W_n}$ 
be the ring homomorphism given by
\begin{equation}
\label{eq:pi}
\tilde{\pi}_n (x_i)
 =
\begin{cases}
 x_i &\text{if $1 \le i \le n$,} \\
 x_{i-n}^{-1} &\text{if $n+1 \le i \le 2n$,} \\
 0 &\text{if $i \ge 2n+1$.}
\end{cases}
\end{equation}
For each strict partition $\lambda$, 
the symmetric functions $\Pfun^C_\lambda$ and $\Qfun^C_\lambda \in \Lambda$ 
are the unique ones satisfying
\begin{equation}
\label{eq:uPC=PC}
\tilde{\pi}_n( \Pfun^C_\lambda )
 = 
P^C_\lambda(x_1, \dots, x_n),
\quad
\tilde{\pi}_n( \Qfun^C_\lambda )
 = 
Q^C_\lambda(x_1, \dots, x_n).
\end{equation}
\end{prop}

Note that $P^C_\lambda(x_1, \dots, x_n) = Q^C_\lambda(x_1, \dots, x_n) = 0$ if $l(\lambda) > n$ 
(Lemma~\ref{lem:PC=0}).

\begin{demo}{Proof}
We have $\tilde{\pi}_n(\Qfun^C_{(r)}) = Q^C_{(r)}(x_1, \dots, x_n)$ by (\ref{eq:GF_q}) and (\ref{eq:GF_length1}).
Then, by (\ref{eq:uQC_length2}) and (\ref{eq:length2}), 
we obtain $\tilde{\pi}_n(\Qfun^C_{(r,s)}) = Q^C_{(r,s)}(x_1, \dots, x_n)$.
Hence we can derive (\ref{eq:uPC=PC}) by comparing (\ref{eq:Schur_QC}) and (\ref{eq:uQC}).
The uniqueness follows from Part (1) of the lemma below.
\end{demo}

\begin{lemma}
\label{lem:proj}
\begin{enumerate}
\item[(1)]
If $f \in \Lambda$ satisfies $\tilde{\pi}_n(f) = 0$ for any large enough $n$, 
then we have $f = 0$.
\item[(2)]
If $f \in \Lambda \otimes \Lambda$ satisfies $(\tilde{\pi}_n \otimes \tilde{\pi}_m)(f) = 0$ 
for any large enough $n$ and $m$, 
then we have $f = 0$.
\end{enumerate}
\end{lemma}

\begin{demo}{Proof}
(1)
Since $\Lambda = \Rat[e_r : r \ge 1]$, 
there exists a positive integer $n$ such that $f \in \Rat[e_1, \dots, e_n]$ and $\tilde{\pi}_n(f) = 0$..
Since $\tilde{\Lambda}^{(n)} = \Rat[x_1^{\pm 1}, \dots, x_n^{\pm 1}]^{W_n}$ 
is a polynomial ring on algebraically independent generators $\tilde{\pi}_n(e_1), \dots, \tilde{\pi}_n(e_n)$, 
we see that the restriction of $\tilde{\pi}_n$ to $\Rat[e_1, \dots, e_n]$ 
gives an isomorphism between $\Rat[e_1, \dots, e_n]$ and $\Rat[x_1^{\pm 1}, \dots, x_n^{\pm 1}]^W$. 
Hence we can conclude $f=0$.

(2) is proved similarly to (1).
\end{demo}

By comparing the definitions of Schur $Q$-functions and universal symplectic $Q$-functions, 
we can show the following proposition.

\begin{prop}
\label{prop:uQC_basis}
\begin{enumerate}
\item[(1)]
For a strict partition $\lambda$, the universal symplectic $P$-/$Q$-functions 
are expressed as linear combinations of Schur $P$-/$Q$-functions 
in the form
\begin{equation}
\label{eq:uQC_by_Q}
\Pfun^C_\lambda
 = 
P_\lambda + \sum_\mu b''_{\lambda,\mu} P_\mu,
\quad
\Qfun^C_\lambda
 = 
Q_\lambda + \sum_\mu b'_{\lambda,\mu} Q_\mu
\end{equation}
respectively, 
where $\mu$ runs over all strict partitions such that $|\mu|<|\lambda|$ and $|\lambda|-|\mu|$ is even.
\item[(2)]
Both $\{ \Pfun^C_\lambda : \lambda \in \SPar \}$ and $\{ \Qfun^C_\lambda : \lambda \in \SPar \}$ 
are bases of $\Gamma$, 
where $\SPar$ is the set of all strict partitions.
\end{enumerate}
\end{prop}

\begin{demo}{Proof}
(1)
Comparing (\ref{eq:Q_length2}) with (\ref{eq:uQC_length2}), 
we have $\Qfun^C_{r,s} - Q_{r,s} \in \bigoplus_{i \ge 1} \Gamma_{r+s-2i}$.
Hence, by expanding the Pfaffians in (\ref{eq:Q}) and (\ref{eq:uQC}), 
we see that 
$\Qfun^C_\lambda - Q_\lambda \in \bigoplus_{i \ge 1} \Gamma_{|\lambda|-2i}$.
Then by using Proposition~\ref{prop:Q_basis}, we can express $\Qfun^C_\lambda$ 
as a linear combination of Schur $Q$-functions in the form (\ref{eq:uQC_by_Q}).

(2) 
The claim follows from Proposition~\ref{prop:Q_basis} and (\ref{eq:uQC_by_Q}).
\end{demo}

\subsection{%
Universal skew symplectic $P$-functions
}

Now we introduce a skew version of universal symplectic $Q$-functions.
For sequences of nonnegative integers $\alpha = (\alpha_1, \dots, \alpha_r)$ 
and $\beta = (\beta_1, \dots, \beta_s)$, we put
$$
\tilde{K}_\alpha
 =
\Big( \Qfun^C_{(\alpha_i,\alpha_j)} \Big)_{1 \le i, j \le r},
\quad
\tilde{M}_{\alpha/\beta}
 =
\Big( \Qfun^C_{(\alpha_i-\beta_{s+1-j})} \Big)_{1 \le i \le r, 1 \le j \le s},
$$
where $\Qfun^C_{(0)} = 1$ and $\Qfun^C_{(k)} = 0$ for $k<0$.
Given two strict partition $\lambda$ and $\mu$, we define 
\begin{equation}
\label{eq:uskewQC}
\Qfun^C_{\lambda/\mu}
 =
\begin{cases}
 \Pf \begin{pmatrix}
  \tilde{S}_\lambda & \tilde{M}_{\lambda/\mu} \\
  -\trans \tilde{M}_{\lambda/\mu} & O
 \end{pmatrix}
&\text{if $l(\lambda)+l(\mu)$ is even,}
\\
 \Pf \begin{pmatrix}
  \tilde{S}_\lambda & \tilde{M}_{\lambda/\mu^0} \\
  -\trans \tilde{M}_{\lambda/\mu^0} & O
 \end{pmatrix}
&\text{if $l(\lambda)+l(\mu)$ is odd.}
\end{cases}
\end{equation}
This definition is similar to the Pragacz--J\'ozefiak--Nimmo formula 
for skew $Q$-functions (see \cite[Theorem~1]{PJ91} and \cite[(2.22)]{Nimmo90}).
And we define
\begin{equation}
\label{eq:uskewPC}
\Pfun^C_{\lambda/\mu}
 =
2^{-l(\lambda)+l(\mu)} \Qfun^C_{\lambda/\mu}.
\end{equation}
We call $\Qfun^C_{\lambda/\mu}$ and $\Pfun^C_{\lambda/\mu}$ 
the universal skew symplectic $Q$-function and $P$-function respectively.
For a finite number of variables $\vectx = (x_1, \dots, x_n)$, we put
$$
Q^C_{\lambda/\mu}(\vectx)
 =
\tilde{\pi}_n ( \Qfun^C_{\lambda/\mu} ),
\quad
P^C_{\lambda/\mu}(\vectx)
 =
\tilde{\pi}_n ( \Pfun^C_{\lambda/\mu} ).
$$
Comparing (\ref{eq:uQC}) with (\ref{eq:uskewQC}), we obtain $\Qfun^C_{\lambda/\emptyset} = \Qfun^C_\lambda$.
Since $\Qfun^C_{(k)} = 0$ for $k<0$, 
we can use the same arguments as in the proof of \cite[Propositions~4.4]{Okada20} 
to prove the following proposition.

\begin{prop}
\label{prop:uskewQC}
For two strict partitions $\lambda$ and $\mu$, 
we have $\Qfun^C_{\lambda/\mu} = 0$ unless $\lambda \supset \mu$, i.e., $S(\lambda) \supset S(\mu)$.
\end{prop}

The following proposition justifies the name ``skew'' symplectic $Q$-functions.

\begin{prop}
\label{prop:sep_var}
Let $X$ and $Y$ be two disjoint sets of infinitely many variables.
For strict partitions $\lambda$ and $\nu$, we have
\begin{gather}
\label{eq:sep_var1}
\Qfun^C_\lambda(X \cup Y)
 =
\sum_\mu \Qfun^C_{\lambda/\mu}(X) \Qfun^C_\mu(Y),
\\
\label{eq:sep_var2}
\Qfun^C_{\lambda/\nu}(X \cup Y)
 =
\sum_\mu \Qfun^C_{\lambda/\mu}(X) \Qfun^C_{\mu/\nu}(Y),
\end{gather}
where $\mu$ runs over all strict partitions.
\end{prop}

For the proof of this proposition, we need the following relations, 
which correspond the cases where $l(\lambda) = 1$ and $2$ of (\ref{eq:sep_var1}).

\begin{lemma}
\label{lem:sep_var}
For nonnegative integers $r$ and $s$, we have
\begin{gather}
\label{eq:sep_var_length1}
\Qfun^C_{(r)}(X \cup Y)
 = 
\sum_{k \ge 0} \Qfun^C_{(k)}(X) \Qfun^C_{(r-k)}(Y),
\\
\label{eq:sep_var_length2}
\Qfun^C_{(r,s)}(X \cup Y)
 =
\Qfun^C_{(r,s)}(X)
 +
\sum_{k, l \ge 0} \Qfun^C_{(r-k)}(X) \Qfun^C_{(s-l)}(Y) \Qfun^C_{(k,l)}(Y).
\end{gather}
\end{lemma}

\begin{demo}{Proof}
Since $\Qfun^C_{(r)} = q_r$ by the definition (\ref{eq:uQC_length1}), 
Equation (\ref{eq:sep_var_length1}) follows from (\ref{eq:GF_q}).
By Lemma~\ref{lem:proj} (2), it is enough to prove (\ref{eq:sep_var_length2}) for finitely many variables.
Let $\vectx = (x_1, \dots, x_n)$ and $\vecty = (y_1, \dots, y_m)$, 
and consider the generating function.
By using Proposition~\ref{prop:length12}, we have
\begin{align*}
&
\sum_{r,s \ge 0}
 \left(
  Q^C_{(r,s)}(\vectx)
  +
  \sum_{k, l \ge 0} Q^C_{(r-k)}(\vectx) Q^C_{(s-l)}(\vectx) Q^C_{(k,l)}(\vecty)
 \right)
 z^r w^s
\\
&
 =
\sum_{r, s \ge 0} Q^C_{(r,s)}(\vectx) z^r w^s
 +
\left( \sum_{a \ge 0} Q^C_{(a)}(\vectx) z^a \right)
\left( \sum_{b \ge 0} Q^C_{(b)}(\vectx) z^b \right)
\left( \sum_{k,l \ge 0} Q^C_{(k,l)}(\vecty) z^k w^l \right)
\\
&
 =
\frac{ (z-w)(1-zw) }
     { (z+w)(1+zw) }
\left( \tilde{\Pi}_z(\vectx) \tilde{\Pi}_w(\vectx) - 1 \right)
+
\tilde{\Pi}_z(\vectx) \tilde{\Pi}_w(\vectx)
\cdot
\frac{ (z-w)(1-zw) }
     { (z+w)(1+zw) }
\left( \tilde{\Pi}_z(\vecty) \tilde{\Pi}_w(\vecty) - 1 \right)
\\
&
=
\frac{ (z-w)(1-zw) }
     { (z+w)(1+zw) }
\left( \tilde{\Pi}_z(\vectx,\vecty) \tilde{\Pi}_w(\vectx,\vecty) - 1 \right)
\\
&
=
\sum_{r,s \ge 0}
 Q^C_{(r,s)}(\vectx,\vecty) z^r w^s.
\end{align*}
Hence we obtain (\ref{eq:sep_var_length2}).
\end{demo}

\begin{demo}{Proof of Proposition~\ref{prop:sep_var}}
First we derive (\ref{eq:sep_var2}) from (\ref{eq:sep_var1}).
If $Z$ is another set of infinitely many variables, then by using (\ref{eq:sep_var1}) we have
\begin{align*}
\sum_\nu \Qfun^C_{\lambda/\nu}(X \cup Y) \Qfun^C_\nu(Z)
 &=
\Qfun^C_\lambda(X \cup Y \cup Z)
 =
\sum_\mu \Qfun^C_{\lambda/\mu}(X) \Qfun^C_\mu(Y \cup Z)
\\
 &=
\sum_{\mu,\nu} \Qfun^C_{\lambda/\mu}(X) \Qfun^C_{\mu/\nu}(Y) \Qfun^C_\nu(Z).
\end{align*}
By equating the coefficients of $\Qfun^C_\nu(Z)$, we obtain (\ref{eq:sep_var2}).

Now we prove (\ref{eq:sep_var1}).
We use the following Pfaffian version of Ishikawa--Wakayama's minor-summation formula 
(\cite[Theorem~3.4]{Okada19}).
Let $r$ be an even integer and $[r] = \{ 1, \dots, r \}$, $\Nat$ the set of nonnegative integers.
Let $A$ be an $r \times r$ skew-symmetric matrix with rows and columns indexed by $[r]$, 
$B$ a skew-symmetric matrix with rows and columns indexed by $\Nat$, 
and $S$ an $r$-rowed matrix with rows indexed by $[r]$ and columns indexed by $\Nat$, 
Then we have
\begin{equation}
\label{eq:Pf-IW}
\sum_I 
 \Pf B(I) 
 \Pf \begin{pmatrix}
  A & S([r];I) \\
  - \trans S([r];I) & O
 \end{pmatrix}
=
\Pf \left(
 A - S B \trans S
\right),
\end{equation}
where $I$ runs over all even-element subsets of $\Nat$. 
Here $B(I)$ stands for the skew-symmetric submatrix of $B$ obtained by taking rows/columns with indices in $I$,
and $S([r];I)$ for the submatrix of $S$ consisting for columns with indices in $I$.

First we consider the case $l = l(\lambda)$ is even.
In this case, we apply the above formula (\ref{eq:Pf-IW}) to the matrices
$$
A = \Big( \Qfun^C_{(\lambda_i,\lambda_j)}(X) \Big)_{1 \le i, j \le l},
\quad
B = \Big( - \Qfun^C_{(i,j)}(Y) \Big)_{i,j \ge 0},
\quad
S = \Big( \Qfun^C_{(\lambda_i-j)}(X) \Big)_{1 \le i \le l, j \ge 0}.
$$
Strict partitions $\mu$ are in bijection with even-element subsets 
$I = \{ \mu_1, \dots, \mu_m \}$ with $m = l(\mu)$ or $l(\mu)+1$.
If $I$ corresponds to $\mu$, then by (\ref{eq:uQC}) and (\ref{eq:uskewQC}) we have
\begin{gather*}
\Pf B(I)
 = 
\Pf \Big( - \Qfun^C_{(\mu_{m+1-i},\mu_{m+1-j})}(Y) \Big)_{1 \le i, j \le m}
 =
\Pf \Big( \Qfun^C_{(\mu_i,\mu_j)}(Y) \Big)_{1 \le i, j \le m}
 =
\Qfun^C_\mu(Y),
\\
\Pf \begin{pmatrix}
 A & S([r];I) \\
 - \trans S([r];I) & O
\end{pmatrix}
 =
\Qfun^C_{\lambda/\mu}(X).
\end{gather*}
By using (\ref{eq:sep_var_length2}), we see that the $(i,j)$-entry of $A - S B \trans S$ is equal to
$$
\Qfun^C_{(\lambda_i,\lambda_j)}(X)
 + 
\sum_{k,l} \Qfun^C_{(\lambda_i-k)}(X) \Qfun^C_{(\lambda_j-l)}(X) \Qfun^C_{(k,l)}(Y)
 =
\Qfun^C_{(\lambda_i,\lambda_j)}(X \cup Y),
$$
hence we obtain $\Pf (A - S B \trans S) = \Qfun^C_\lambda(X \cup Y)$ by (\ref{eq:uQC}).

Next we consider the case where $l = l(\lambda)$ is odd.
In this case, we apply the above formula (\ref{eq:Pf-IW}) to the matrices
$$
A = \Big( \Qfun^C_{(\lambda_i,\lambda_j)}(X) \Big)_{1 \le i, j \le l+1},
\quad
B = \Big( - \Qfun^C_{(i,j)}(Y) \Big)_{i,j \ge 0},
\quad
S = \Big( \Qfun^C_{(\lambda_i-j)}(X) \Big)_{1 \le i \le l+1, j \ge 0}.
$$
If an even-element subset $I \subset \Nat$ corresponds to a strict partition $\mu$, then we have
$$
\Pf B(I)
 =
\Qfun^C_\mu(Y),
\quad
\Pf \begin{pmatrix}
 A & S([r];I) \\
 - \trans S([r];I) & O
\end{pmatrix}
 =
\Qfun^C_{\lambda/\mu}(X).
$$
By using (\ref{eq:sep_var_length2}), 
we see that the $(i,j)$-entry ($1 \le i < j \le l$) of $A - S B \trans S$ is equal to 
$\Qfun^C_{(\lambda_i,\lambda_j)}(X \cup Y)$.
By using $\Qfun^C_{(\lambda_{l+1}-j)} = \delta_{j,0}$, 
$\Qfun^C_{(r,0)} = \Qfun^C_{(r)}$ ($r>0$), $\Qfun^C_{(0,0)} = 0$ 
and (\ref{eq:sep_var_length1}), 
we see that the $(i,l+1)$ entry of $A - S B \trans S$ is equal to
$$
\Qfun^C_{(\lambda_i,0)}(X) + \sum_{k \ge 1} \Qfun^C_{(\lambda_i-k)}(X) \Qfun^C_{(k,0)}(Y)
 =
\Qfun^C_{(\lambda_i)}(X \cup Y)
 =
\Qfun^C_{(\lambda_i,0)}(X \cup Y).
$$
Hence we have $\Pf (A - \trans S B S) = \Qfun^C_\lambda(X \cup Y)$ by (\ref{eq:uQC}).
\end{demo}

\section{%
Tableau description
}

In this section, we give a proof of a tableau description of symplectic $Q$-functions, 
which was originally conjectured by King--Hamel \cite{KH07}.
One of the keys to the proof is a determinant formula for skew symplectic $Q$-functions 
(see Lemma~\ref{lem:algQ} (3)).

Let $\lambda$ and $\mu$ be strict partitions, 
and $S(\lambda)$ and $S(\mu)$ the corresponding shifted diagrams respectively. 
If $S(\lambda) \supset S(\mu)$, then we write $\lambda \supset \mu$, 
and define the \emph{skew shifted diagram} $S(\lambda/\mu)$ as 
the set-theoretical difference $S(\lambda/\mu) = S(\lambda) \setminus S(\mu)$.

\begin{definition}
\label{def:tableaux}
(Hamel--King \cite{HK07})
Let $\lambda$ and $\mu$ be strict partitions such that $\lambda \supset \mu$.
A \emph{symplectic primed shifted tableaux of shape $\lambda/\mu$} 
is a filling of the cells of $S(\lambda/\mu)$ with entries from the totally ordered set
$$
\mathcal{A}_n
=
\{
1' < 1 < \overline{1}' < \overline{1} <
2' < 2 < \overline{2}' < \overline{2} <
\dots <
n' < n < \overline{n}' < \overline{n}
\}
$$
satisfying the following 5 conditions:
\begin{enumerate}
\item[(T1)]
the entries in each row weakly increase from left to right;
\item[(T2)]
the entries in each column weakly increase from top to bottom;
\item[(T3)]
each row contains at most one $k'$ and at most one $\overline{k}'$;
\item[(T4)]
each column contains at most one $k$ and at most one $\overline{k}$;
\item[(T5)]
the main diagonal contains at most one entry from $\{ k', k, \overline{k}', \overline{k} \}$.
\end{enumerate}
For a symplectic primed shifted tableaux $T$ of shape $\lambda/\mu$, we define its weight $\vectx^T$ 
by putting
$$
\vectx^T
 =
\prod_{k=1}^n x_i^{m(k')+m(k)-m(\overline{k}')-m(\overline{k})},
$$
where $m(\gamma)$ is the multiplicity of $\gamma$ in $T$.
We denote by $\QTab^C_n(\lambda/\mu)$ the set of all symplectic primed shifted tableaux 
of shape $\lambda/\mu$ with entries from $\mathcal{A}_n$.
We write $\QTab^C_n(\lambda)$ for $\QTab^C_n(\lambda/\emptyset)$, 
where $\emptyset$ is the unique partition of $0$.
\end{definition}

This definition is a symplectic analogue of primed shifted tableaux for Schur $Q$-functions 
(see \cite[(8.16')]{Macdonald95})
and a $Q$-function analogue of symplectic tableaux for symplectic Schur functions 
(see \cite[Section~4]{King76}).

For example,
$$
T =
\raisebox{-25pt}{
\setlength{\unitlength}{1.4pt}
\begin{picture}(80,40)
\put(0,40){\line(1,0){80}}
\put(0,30){\line(1,0){80}}
\put(10,20){\line(1,0){60}}
\put(20,10){\line(1,0){50}}
\put(30,0){\line(1,0){20}}
\put(0,40){\line(0,-1){10}}
\put(10,40){\line(0,-1){20}}
\put(20,40){\line(0,-1){30}}
\put(30,40){\line(0,-1){40}}
\put(40,40){\line(0,-1){40}}
\put(50,40){\line(0,-1){40}}
\put(60,40){\line(0,-1){30}}
\put(70,40){\line(0,-1){30}}
\put(80,40){\line(0,-1){10}}
\put(0,30){\makebox(10,10){$1$}}
\put(10,30){\makebox(10,10){$1$}}
\put(20,30){\makebox(10,10){$2'$}}
\put(30,30){\makebox(10,10){$2$}}
\put(40,30){\makebox(10,10){$3$}}
\put(50,30){\makebox(10,10){$\overline{3}$}}
\put(60,30){\makebox(10,10){$\overline{3}$}}
\put(70,30){\makebox(10,10){$5$}}
\put(10,20){\makebox(10,10){$2'$}}
\put(20,20){\makebox(10,10){$\overline{2}'$}}
\put(30,20){\makebox(10,10){$3'$}}
\put(40,20){\makebox(10,10){$\overline{3}'$}}
\put(50,20){\makebox(10,10){$4'$}}
\put(60,20){\makebox(10,10){$4$}}
\put(20,10){\makebox(10,10){$3'$}}
\put(30,10){\makebox(10,10){$3$}}
\put(40,10){\makebox(10,10){$\overline{3}'$}}
\put(50,10){\makebox(10,10){$4'$}}
\put(60,10){\makebox(10,10){$\overline{4}$}}
\put(30,0){\makebox(10,10){$\overline{5}'$}}
\put(40,0){\makebox(10,10){$\overline{5}$}}
\end{picture}
}
$$
is a symplectic primed shifted tableau of shape $(8,6,5,2)$ with weight 
$\vectx^T = x_1^2 x_2^2 x_4^2 x_5^{-1}$.

The aim of this section is to prove the following theorem, which was conjectured by 
King and Hamel \cite[Conjecture~3.1]{KH07}.

\begin{theorem}
\label{thm:tableaux}
Let $\vectx = (x_1, \dots, x_n)$ be a sequence of $n$ indeterminates.
For strict partitions $\lambda$ and $\mu$, we have
\begin{align}
\label{eq:tableaux}
Q^C_{\lambda/\mu}(\vectx)
 &=
\sum_{T \in \QTab^C_n(\lambda/\mu)} \vectx^T,
\\
\label{eq:tableaux_P}
P^C_{\lambda/\mu}(\vectx)
 &=
\sum_{T \in \PTab^C_n(\lambda/\mu)} \vectx^T,
\end{align}
where $\PTab^C_n(\lambda/\mu)$ is the subset of $\QTab^C_n(\lambda/\mu)$ 
consisting of $T \in \QTab^C_n(\lambda/\mu)$ with no primed letter on the main diagonal.
\end{theorem}

In order to prove this theorem, we put
$$
Q^{\text{tab}}_{\lambda/\mu}(\vectx)
 =
\sum_{T \in \QTab^C_n(\lambda/\mu)} \vectx^T
$$
for two strict partitions $\lambda \supset \mu$, 
and show that the generating functions $Q^{\text{tab}}_{\lambda/\mu}(\vectx)$ satisfy the same relations 
as symplectic $Q$-functions $Q^C_{\lambda/\mu}(\vectx) = \tilde{\pi}_n( \Qfun^C_{\lambda/\mu} )$.

\begin{lemma}
\label{lem:algQ}
Let $\lambda$ and $\mu$ be strict partitions such that $\lambda \supset \mu$.
\begin{enumerate}
\item[(1)]
We have
$$
Q^C_{\lambda/\nu}(x_1, \dots, x_n)
 =
\sum \prod_{i=1}^n Q^C_{\mu^{(i)}/\mu^{(i-1)}} (x_i),
$$
where the sum is taken over all sequences $\mu = \mu^{(0)} \subset \mu^{(1)} \subset \dots \subset 
\mu^{(n-1)} \subset \mu^{(n)} = \lambda$.
\item[(2)]
For a single variable $x_1$, we have $Q^C_{\lambda/\mu}(x_1) = 0$ unless $l(\lambda) - l(\mu) \le 1$.
\item[(3)]
If $l(\lambda) - l(\mu) \le 1$, then we have
$$
Q^C_{\lambda/\mu}(x_1)
 =
\det \Big( Q^C_{(\lambda_i - \mu_j)}(x_1) \Big)_{1 \le i, j \le l(\lambda)},
$$
where $Q^C_{(r)}(x_1) = 0$ for $r<0$.
\end{enumerate}
\end{lemma}

\begin{demo}{Proof}
(1) is obtained by iteratively applying (\ref{eq:sep_var1}). 
Since $Q^C_{(r,s)}(x_1) = 0$, 
(2) and (3) follows from the definition (\ref{eq:uskewQC}) and (\ref{eq:Pf-Laplace}).
\end{demo}

\begin{lemma}
\label{lem:combQ}
Let $\lambda$ and $\mu$ be strict partitions such that $\lambda \supset \mu$.
\begin{enumerate}
\item[(1)]
We have
$$
Q^{\text{tab}}_{\lambda/\nu}(x_1, \dots, x_n)
 =
\sum \prod_{i=1}^n Q^{\text{tab}}_{\mu^{(i)}/\mu^{(i-1)}} (x_i),
$$
where the sum is taken over all sequences $\mu = \mu^{(0)} \subset \mu^{(1)} \subset \dots \subset 
\mu^{(n-1)} \subset \mu^{(n)} = \lambda$.
\item[(2)]
For a single variable $x_1$, we have $Q^{\text{tab}}_{\lambda/\mu}(x_1) = 0$ unless $l(\lambda) - l(\mu) \le 1$.
\item[(3)]
If $l(\lambda) - l(\mu) \le 1$, then we have
$$
Q^{\text{tab}}_{\lambda/\mu}(x_1)
 =
\det \Big( Q^{\text{tab}}_{(\lambda_i - \mu_j)}(x_1) \Big)_{1 \le i, j \le l(\lambda)},
$$
where $Q^{\text{tab}}_{(r)}(x) = 0$ for $r<0$.
\end{enumerate}
\end{lemma}

\begin{demo}{Proof}
(1)
The claim is obtained by decomposing a primed shifted tableaux into 
subtableaux consisting of $i'$, $i$, $\overline{i}'$ and $\overline{i}$ for $i=1, \dots, n$, 

(2)
Condition (T5) in Definition~\ref{def:tableaux} implies that, 
if $l(\lambda) - l(\mu) \ge 2$, then there are no primed shifted tableaux of shape $\lambda/\mu$ 
with entries from $\{ 1', 1, \overline{1}', \overline{1} \}$.

(3)
We appeal to the Lindstr\"om--Gessel--Viennot lemma together with 
a bijection between symplectic primed shifted tableaux 
and non-intersecting lattice paths.
Let $G = (V,E)$ the directed graph with vertex set
$$
V = \{ A_i = (i,0) : i \ge 0 \} \cup \{ B_i = (i,1) : i \ge 0 \} \cup \{ C_i = (i,2) : i \ge 0 \},
$$
and directed edge set
\begin{multline*}
E
 = 
\{ (A_i, B_i), (B_i, C_i) : i \ge 0 \}
\cup
\{ (A_i, B_{i+1}), (B_i, C_{i+1}) : i \ge 0 \}
\\
\cup
\{ (B_i, B_{i+1}), (C_i, C_{i+1}) : i \ge 0 \}.
\end{multline*}
See Figure~\ref{fig:G}.
\begin{figure}
\centering
\setlength\unitlength{1.5pt}
\begin{picture}(110,60)(-10,-10)
\put(0,40){\circle*{2}}
\put(20,40){\circle*{2}}
\put(40,40){\circle*{2}}
\put(60,40){\circle*{2}}
\put(80,40){\circle*{2}}
\put(0,20){\circle*{2}}
\put(20,20){\circle*{2}}
\put(40,20){\circle*{2}}
\put(60,20){\circle*{2}}
\put(80,20){\circle*{2}}
\put(0,0){\circle*{2}}
\put(20,0){\circle*{2}}
\put(40,0){\circle*{2}}
\put(60,0){\circle*{2}}
\put(80,0){\circle*{2}}
\put(85,35){\makebox(10,10){$\cdots$}}
\put(85,15){\makebox(10,10){$\cdots$}}
\put(85,-5){\makebox(10,10){$\cdots$}}
\put(0,40){\vector(1,0){19}}
\put(20,40){\vector(1,0){19}}
\put(40,40){\vector(1,0){19}}
\put(60,40){\vector(1,0){19}}
\put(0,20){\vector(0,1){19}}
\put(20,20){\vector(0,1){19}}
\put(40,20){\vector(0,1){19}}
\put(60,20){\vector(0,1){19}}
\put(80,20){\vector(0,1){19}}
\put(0,20){\vector(1,1){19}}
\put(20,20){\vector(1,1){19}}
\put(40,20){\vector(1,1){19}}
\put(60,20){\vector(1,1){19}}
\put(0,20){\vector(1,0){19}}
\put(20,20){\vector(1,0){19}}
\put(40,20){\vector(1,0){19}}
\put(60,20){\vector(1,0){19}}
\put(0,0){\vector(0,1){19}}
\put(20,0){\vector(0,1){19}}
\put(40,0){\vector(0,1){19}}
\put(60,0){\vector(0,1){19}}
\put(80,0){\vector(0,1){19}}
\put(0,0){\vector(1,1){19}}
\put(20,0){\vector(1,1){19}}
\put(40,0){\vector(1,1){19}}
\put(60,0){\vector(1,1){19}}
\put(-5,40){\makebox(10,10){$C_0$}}
\put(15,40){\makebox(10,10){$C_1$}}
\put(35,40){\makebox(10,10){$C_2$}}
\put(55,40){\makebox(10,10){$C_3$}}
\put(75,40){\makebox(10,10){$C_4$}}
\put(-10,20){\makebox(10,10){$B_0$}}
\put(10,20){\makebox(10,10){$B_1$}}
\put(30,20){\makebox(10,10){$B_2$}}
\put(50,20){\makebox(10,10){$B_3$}}
\put(70,20){\makebox(10,10){$B_4$}}
\put(-5,-10){\makebox(10,10){$A_0$}}
\put(15,-10){\makebox(10,10){$A_1$}}
\put(35,-10){\makebox(10,10){$A_2$}}
\put(55,-10){\makebox(10,10){$A_3$}}
\put(75,-10){\makebox(10,10){$A_4$}}
\end{picture}
\caption{Graph $G$}
\label{fig:G}
\end{figure}
For two vertices $A_s$, $C_r \in V$, a \emph{lattice path} from $A_s$ to $C_r$ 
is a sequence $(S_0, S_1, \dots, S_m)$ 
of vertices of $G$ such that $S_0 = A_s$, $S_m = C_r$ and $(S_k, S_{k+1}) \in E$ for $0 \le k \le m-1$.
We denote by $\LL(s;r)$ the set of all lattice paths from $A_s$ to $C_r$.
A family $(P_1, \dots, P_l)$ of lattice paths in $G$ is called 
\emph{non-intersecting} if no two of them have a vertex in common.
For two strict partitions $\lambda$ and $\mu$ such that $l(\lambda) = l$ and $l(\mu) = l$ or $l-1$, 
we denote by $\LL_0(\mu;\lambda)$ the set of all non-intersecting lattice paths 
$(P_1, \dots, P_l)$ with $P_i \in \LL(\mu_i;\lambda_i)$.
Note that $\mu_l = 0$ if $l(\mu) = l(\lambda)-1$.

We consider a weight function $\wt : E \to \Int[x_1^{\pm 1}]$ given by
\begin{gather*}
\wt(A_i,B_i) = \wt(B_i,C_i) = 1,
\\
\wt(A_i,B_{i+1}) = \wt(B_i,B_{i+1}) = x_1,
\\
\wt(B_i,C_{i+1}) = \wt(C_i,C_{i+1}) = x_1^{-1}.
\end{gather*}
Then we define the weight $\wt(P)$ of a lattice path $P$ to be 
the product of the weights of its individual steps, 
and put $\wt(P_1, \dots, P_l) = \prod_{i=1}^l \wt(P_i)$.

There exists a weight-preserving bijection from $\QTab^C_1(\lambda/\mu)$ to $\LL_0(\mu;\lambda)$.
If a tableau $T \in \QTab^C_1(\lambda/\mu)$ corresponds to a family of lattice paths $(P_1, \dots, P_l)$, 
then the entries of the $i$th row of $T$ are obtained by reading 
the diagonal and horizontal edges of the $i$th lattice path $P_i$ 
from left to right, 
and assigning $1'$, $1$, $\overline{1}'$ and $\overline{1}$ 
to the edges $(A_k, B_{k+1})$, $(B_k,B_{k+1})$, $(B_k,C_{k+1})$ and $(C_k,C_{k+1})$ respectively.
See Figure~\ref{fig:bijection} for an example of the correspondence 
in the case $\lambda = (6,5,2)$ and $\mu = (3,2)$.
\begin{figure}
\centering
\setlength\unitlength{1.5pt}
\raisebox{22.5pt}{
\begin{picture}(60,30)
\dashline{1.5}(0,30)(30,30)
\put(30,30){\line(1,0){30}}
\dashline{1.5}(0,20)(30,20)
\put(30,20){\line(1,0){30}}
\dashline{1.5}(10,10)(20,10)
\put(20,10){\line(1,0){40}}
\put(20,0){\line(1,0){20}}
\dashline{1.5}(0,30)(0,20)
\dashline{1.5}(10,30)(10,10)
\dashline{1.5}(20,30)(20,10)
\put(20,10){\line(0,-1){10}}
\put(30,30){\line(0,-1){30}}
\put(40,30){\line(0,-1){30}}
\put(50,30){\line(0,-1){20}}
\put(60,30){\line(0,-1){20}}
\put(30,20){\makebox(10,10){$1'$}}
\put(40,20){\makebox(10,10){$1$}}
\put(50,20){\makebox(10,10){$1$}}
\put(30,10){\makebox(10,10){$1$}}
\put(40,10){\makebox(10,10){$\overline{1}'$}}
\put(50,10){\makebox(10,10){$\overline{1}$}}
\put(20,0){\makebox(10,10){$1$}}
\put(30,0){\makebox(10,10){$\overline{1}$}}
\end{picture}
}
\quad
\raisebox{37.5pt}{
\begin{picture}(20,10)
\put(0,0){\makebox(20,10){$\longleftrightarrow$}}
\end{picture}
}
\quad
\begin{picture}(140,60)(-10,-10)
\put(0,40){\circle*{2}}
\put(20,40){\circle*{2}}
\put(40,40){\circle*{2}}
\put(60,40){\circle*{2}}
\put(80,40){\circle*{2}}
\put(100,40){\circle*{2}}
\put(120,40){\circle*{2}}
\put(0,20){\circle*{2}}
\put(20,20){\circle*{2}}
\put(40,20){\circle*{2}}
\put(60,20){\circle*{2}}
\put(80,20){\circle*{2}}
\put(100,20){\circle*{2}}
\put(120,20){\circle*{2}}
\put(0,0){\circle*{2}}
\put(20,0){\circle*{2}}
\put(40,0){\circle*{2}}
\put(60,0){\circle*{2}}
\put(80,0){\circle*{2}}
\put(100,0){\circle*{2}}
\put(120,0){\circle*{2}}
\put(0,0){\vector(0,1){19}}
\put(0,20){\vector(1,0){19}}
\put(20,20){\vector(0,1){19}}
\put(20,40){\vector(1,0){19}}
\put(100,20){\vector(1,0){19}}
\put(80,40){\vector(1,0){19}}
\put(40,0){\vector(0,1){19}}
\put(40,20){\vector(1,0){19}}
\put(60,20){\vector(1,1){19}}
\put(80,40){\vector(1,0){19}}
\put(60,0){\vector(1,1){19}}
\put(80,20){\vector(1,0){19}}
\put(100,20){\vector(1,0){19}}
\put(120,20){\vector(0,1){19}}
\put(35,40){\makebox(10,10){$C_2$}}
\put(95,40){\makebox(10,10){$C_5$}}
\put(115,40){\makebox(10,10){$C_6$}}
\put(-5,-10){\makebox(10,10){$A_0$}}
\put(35,-10){\makebox(10,10){$A_2$}}
\put(55,-10){\makebox(10,10){$A_3$}}
\end{picture}
\caption{Bijection}
\label{fig:bijection}
\end{figure}
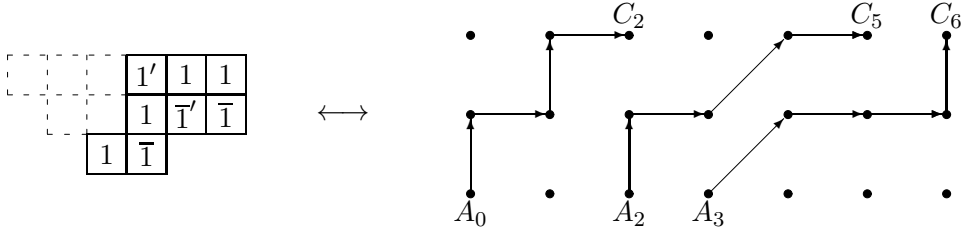

By using this bijection and the Lindstr\"om--Gessel--Viennot Lemma, we have
\begin{align*}
Q^{\text{tab}}_{\lambda/\mu}(x_1) 
 &= 
\sum_{T \in \QTab^C_1(\lambda/\mu)} x_1^{m(1')+m(1)-m(\overline{1}')-m(\overline{1})}
 =
\sum_{(P_1, \dots, P_l) \in \LL_0(\mu;\lambda)} \prod_{i=1}^l \wt(P_i)
\\
 &=
\det \left(
 \sum_{P \in \LL(\mu_j;\lambda_i)} \wt(P)
\right)_{1 \le i, j \le l}
 =
\det \Big(
 Q^{\text{tab}}_{(\lambda_i-\mu_j)} (x_1)
\Big)_{1 \le i, j \le l}.
\end{align*}
\end{demo}

Now we are ready to complete the proof of Theorem~\ref{thm:tableaux}.

\begin{demo}{Proof of Theorem~\ref{thm:tableaux}}
By Lemmas~\ref{lem:algQ} and \ref{lem:combQ}, it is enough to show 
$Q^C_{(r)}(x_1) = Q^{\text{tab}}_{(r)}(x_1)$.

Symplectic primed shifted tableaux $T \in \QTab^C_1((r))$ of shape $(r)$ are in bijection with 
quadruples $(a,b,c,d) = (m(1'), m(1), m(\overline{1}'), m(\overline{1}))$ 
of nonnegative integers such that $a$, $c \in \{ 0, 1 \}$
and $a+b+c+d = r$.
Hence the generating function of $Q^{\text{tab}}_{(r)}(x_1)$ is given by
$$
\sum_{r \ge 0} Q^{\text{tab}}_{(r)}(x_1) z^r
 =
\sum_{a=0}^1 \sum_{b=0}^\infty \sum_{c=0}^1 \sum_{d=0}^\infty
 x_1^{a+b-c-d} z^{a+b+c+d}
 =
\frac{ 1 + x_1 z }{ 1 - x_1 z }
\frac{ 1 + x_1^{-1} z }{ 1 - x_1^{-1} z }.
$$
Comparing the generating function of $Q^C_{(r)}(x_1)$ given by (\ref{eq:GF_length1}), 
we obtain $Q^C_{(r)}(x_1) = Q^{\text{tab}}_{(r)}(x_1)$.
This completes the proof of (\ref{eq:tableaux}).

Since each entry on the main diagonal of $T \in \QTab^C_n(\lambda/\mu)$ 
may be either unprimed or primed, we can obtain (\ref{eq:tableaux_P}) 
from (\ref{eq:tableaux}).
\end{demo}

We conclude this section with a flip symmetry of skew symplectic $Q$-functions.
Let $\delta_r = (r, r-1, \dots, 2, 1)$ be the staircase partition of length $r$.
For a strict partition $\lambda \subset \delta_r$, 
let $\lambda^* = (\lambda^*_1, \lambda^*_2, \dots)$ be the strict partition such that
$\{ \lambda_1, \lambda_2, \dots, \} \sqcup \{ \lambda^*_1, \lambda^*_2, \dots \} 
= \{ 1, 2, \dots, r \}$, 
call it the complement of $\lambda$ in $\delta_r$.
For example, if $\lambda = (5,2) \subset \delta_5 = (5,4,3,2,1)$, 
then we have $\lambda^* = (4,3,1)$.

\begin{prop}
\label{prop:flip}
Let $\lambda$ and $\mu$ be two strict partitions such that $\lambda \supset \mu$.
Let $r$ be a positive integer such that $\delta_r \supset \lambda \supset \mu$, 
and $\lambda^*$ and $\mu^*$ the complements of $\lambda$ and $\mu$ in $\delta_r$ respectively.
Then we have
$$
\Qfun^C_{\lambda/\mu} = \Qfun^C_{\mu^*/\lambda^*}
$$
in $\Lambda$.
\end{prop}

\begin{demo}{Proof}
By Lemma~\ref{lem:proj} (1), it is enough to show 
$Q^C_{\lambda/\mu}(\vectx) = Q^C_{\mu^*/\lambda^*}(\vectx)$ for $\vectx = (x_1, \dots, x_n)$.

Note that the skew shifted diagram $S(\mu^*/\lambda^*)$ is obtained from $S(\lambda/\mu)$ 
by flipping along the anti-diagonal of $S(\delta_r)$.
By replacing $k'$, $k$, $\overline{k}'$ and $\overline{k}$ with $\overline{l}$, $\overline{l}'$, 
$l$ and $l'$ respectively, where $l = n+1-k$, for $k=1, 2, \dots, n$, 
and then flipping the resulting tableau along the anti-diagonal, 
we obtain a bijection $\Phi : \QTab^C_n(\lambda/\mu) \to \QTab^C_n(\mu^*/\lambda^*)$.
For example, if $n=6$, $\lambda = \delta_5 = (5,4,3,2,1)$ and $\mu = (5,2)$, then we have
$$
\setlength{\unitlength}{1.4pt}
T =
\raisebox{-30pt}{
\begin{picture}(50,50)
\dashline{1.5}(0,50)(50,50)
\dashline{1.5}(0,40)(50,40)
\dashline{1.5}(10,30)(30,30)
\dashline{1.5}(0,50)(0,40)
\dashline{1.5}(10,50)(10,30)
\dashline{1.5}(20,50)(20,30)
\dashline{1.5}(30,50)(30,30)
\dashline{1.5}(40,50)(40,40)
\dashline{1.5}(50,50)(50,40)
\put(30,40){\line(1,0){20}}
\put(20,30){\line(1,0){30}}
\put(20,20){\line(1,0){30}}
\put(30,10){\line(1,0){20}}
\put(40,0){\line(1,0){10}}
\put(20,30){\line(0,-1){10}}
\put(30,40){\line(0,-1){30}}
\put(40,40){\line(0,-1){40}}
\put(50,40){\line(0,-1){40}}
\put(30,30){\makebox(10,10){$2'$}}
\put(40,30){\makebox(10,10){$2$}}
\put(20,20){\makebox(10,10){$1$}}
\put(30,20){\makebox(10,10){$2'$}}
\put(40,20){\makebox(10,10){$\overline{3}'$}}
\put(30,10){\makebox(10,10){$\overline{4}$}}
\put(40,10){\makebox(10,10){$\overline{4}$}}
\put(40,0){\makebox(10,10){$5'$}}
\end{picture}
}
\longmapsto
\Phi(T)
=
\raisebox{-30pt}{
\begin{picture}(50,50)
\put(0,50){\line(1,0){40}}
\put(0,40){\line(1,0){40}}
\put(10,30){\line(1,0){30}}
\put(20,20){\line(1,0){10}}
\put(0,50){\line(0,-1){10}}
\put(10,50){\line(0,-1){20}}
\put(20,50){\line(0,-1){30}}
\put(30,50){\line(0,-1){30}}
\put(40,50){\line(0,-1){20}}
\dashline{1.5}(40,50)(50,50)
\dashline{1.5}(40,40)(50,40)
\dashline{1.5}(40,30)(50,30)
\dashline{1.5}(30,20)(50,20)
\dashline{1.5}(30,10)(50,10)
\dashline{1.5}(40,0)(50,0)
\dashline{1.5}(30,20)(30,10)
\dashline{1.5}(40,30)(40,0)
\dashline{1.5}(50,50)(50,0)
\put(0,40){\makebox(10,10){$\overline{2}$}}
\put(10,40){\makebox(10,10){$3'$}}
\put(20,40){\makebox(10,10){$4$}}
\put(30,40){\makebox(10,10){$\overline{5}'$}}
\put(10,30){\makebox(10,10){$3'$}}
\put(20,30){\makebox(10,10){$\overline{5}$}}
\put(30,30){\makebox(10,10){$\overline{5}$}}
\put(20,20){\makebox(10,10){$\overline{6}'$}}
\end{picture}
}.
$$
Let $\phi$ be the ring automorphism of $\Rat[x_1^{\pm 1}, \dots, x_n^{\pm 1}]$ defined by 
$\phi(x_i) = x_{n+1-i}^{-1}$ for $1 \le i \le n$.
Then we have $\vectx^{\Phi(T)} = \phi(\vectx^T)$ for $T \in \QTab^C_n(\lambda/\mu)$.
Since $\phi$ is the action of an element in the Weyl group $W_n$ and $Q^C_{\lambda/\mu}(\vectx)$ is $W_n$-invariant, we have
$$
Q^C_{\lambda/\mu}(\vectx)
 =
\sum_{T \in \QTab^C_n(\lambda/\mu)} \phi(\vectx^T)
 =
\sum_{T \in \QTab^C_n(\lambda/\mu)} \vectx^{\Phi(T)}
 =
\sum_{S \in \QTab^C_n(\mu^*/\lambda^*)} \vectx^S
 =
Q^C_{\mu^*/\lambda^*}(\vectx).
$$
\end{demo}
\section{%
Pieri-type rule
}

This section is devoted to proving a Pieri-type rule, 
which describes the expansion of the product of an arbitrary symplectic $P$-function 
with the symplectic $P$-function corresponding to a one-row partition.

To state the Pieri-type rule for symplectic $P$-functions, we introduce some notation.
For two strict partitions $\lambda$ and $\mu$, we write $\lambda \succ \mu$ 
if $\lambda_1 \ge \mu_1 \ge \lambda_2 \ge \mu_2 \ge \cdots$.
For such a pair of strict partitions, let $a(\lambda,\mu)$ denote the number of connected components 
of the skew shifted diagram $S(\lambda/\mu)$.
Then it is not difficult to see that
\begin{equation}
\label{eq:a}
a(\lambda,\mu)
 =
\# \{ i : 1 \le i \le l-1, \, \lambda_i > \mu_i > \lambda_{i+1} \}
 +
\begin{cases}
 1 &\text{if $\lambda_l > \mu_l$,} \\
 0 &\text{if $\lambda_l = \mu_l$,}
\end{cases}
\end{equation}
where $l$ is the length of $\lambda$.

\begin{theorem}
\label{thm:Pieri}
For strict partitions $\lambda$, $\mu$ and a positive integer $r$, 
we define the Pieri coefficient $\tilde{f}^\lambda_{\mu,(r)}$ 
by the relation
\begin{equation}
\label{eq:Pieri-expansion}
\Pfun^C_\mu \cdot \Pfun^C_{(r)}
 =
\sum_\lambda \tilde{f}^\lambda_{\mu,(r)} \Pfun^C_\lambda.
\end{equation}
Then we have
\begin{equation}
\label{eq:Pieri}
\tilde{f}^\lambda_{\mu,(r)}
 = 
\begin{cases}
\sum_\kappa 2^{a(\mu,\kappa) + a(\lambda,\kappa) - \chi[l(\mu)>l(\kappa)] - 1}
 &\text{if $l(\lambda) = l(\mu)$ or $l(\mu)+1$,} \\
0 &\text{otherwise,}
\end{cases}
\end{equation}
where $\kappa$ runs over all strict partitions satisfying 
$\mu \succ \kappa$, $\lambda \succ \kappa$ and $(|\mu| - |\kappa|) + (|\lambda| - |\kappa|) = r$, 
and
$$
\chi[l(\mu)>l(\kappa)] = \begin{cases}
 1 &\text{if $l(\mu) > l(\kappa)$,} \\
 0 &\text{otherwise.}
\end{cases}
$$
\end{theorem}

This theorem is a symplectic analogue of Morris' Pieri-type rule for Schur $P$-functions 
(\cite[Theorem~1]{Morris64}) 
and a $P$-function analogue of Sundaram's Pieri-type rule for symplectic Schur functions 
(\cite[Theorem~4.1]{Sundaram90}).

Before giving the proof of Theorem~\ref{thm:Pieri}, we present a corollary and an example.
If $r=1$, then we have the following multiplicity-free expansion.

\begin{corollary}
\label{cor:Pieri}
For a strict partition $\mu$, we have
$$
\Pfun^C_\mu \cdot \Pfun^C_{(1)}
 =
\sum_\lambda \Pfun^C_\lambda,
$$
where $\lambda$ runs over all strict partitions 
satisfying one of the following conditions:
\begin{enumerate}
\item[(i)]
$S(\lambda)$ is obtained from $S(\mu)$ by adding one box;
\item[(ii)]
$S(\lambda)$ is obtained from $S(\mu)$ by removing one box, and $l(\lambda) = l(\mu)$.
\end{enumerate}
\end{corollary}

\begin{example}
If $\mu = (4,3,1)$ and $r = 2$, then we have
$$
\Pfun^C_{(4,3,1)} \cdot \Pfun^C_{(2)}
 =
\Pfun^C_{(6, 3, 1)}
 + \Pfun^C_{(5, 4, 1)}
 + 2 \Pfun^C_{(5, 3, 2)}
 + 2 \Pfun^C_{(5, 2, 1)}
 + 3 \Pfun^C_{(4, 3, 1)}
 + \Pfun^C_{(3, 2, 1)}.
$$
If $\lambda = (4,3,1)$, then 
the strict partitions $\kappa$ satisfying $\mu \succ \kappa$, $\lambda \succ \kappa$ 
and $(|\lambda|-|\kappa|) + (|\mu|-|\kappa|) = 2$ 
are $\kappa = (4,2,1)$ and $(4,3)$.
By definition, we have
\begin{gather*}
a((4,3,1),(4,2,1)) = 1, \quad a((4,3,1),(4,3)) = 1,
\\
\chi[l(4,3,1)>l(4,2,1)] = 0, \quad \chi[l(4,3,1)>l(4,3)] = 1,
\end{gather*}
hence we have
$$
\tilde{f}^{(4,3,1)}_{(4,3,1),(2)}
 =
2^{1+1-0-1} + 2^{1+1-1-1}
 =
3.
$$
\end{example}

Now we start the proof of Theorem~\ref{thm:Pieri}.
We take a positive integer $n$ such that $l(\mu) \le n$, 
and apply $\tilde{\pi}_n$ to the both sides of (\ref{eq:Pieri-expansion}).
Then we have
$$
P^C_\mu(\vectx) P^C_{(r)}(\vectx)
 =
\sum_\lambda \tilde{f}^\lambda_{\mu,(r)} P^C_\lambda(\vectx),
$$
where the sum is taken over all strict partitions $\lambda$ of length $\le n$.
By (\ref{eq:GF_length1}), we have
$$
1 + 2 \sum_{r \ge 1} P^C_{(r)}(\vectx) z^r
 =
\tilde{\Pi}_z(\vectx)
 =
\prod_{i=1}^n
 \frac{ (1 + x_i z) (1 + x_i^{-1} z) }
      { (1 - x_i z) (1 - x_i^{-1} z) }.
$$
So we consider the following generating function for the Pieri coefficients 
$\tilde{f}^\lambda_{\mu,(r)}$:
\begin{equation}
\label{eq:GF_Pieri}
u^\lambda_\mu(z)
 =
\delta_{\lambda,\mu} + 2 \sum_{r \ge 1} \tilde{f}^\lambda_{\mu,(r)} z^r.
\end{equation}
Then we have
\begin{equation}
\label{eq:Pieri-expansion2}
P^C_\mu(\vectx) \tilde{\Pi}_z(\vectx)
 =
\sum_\lambda
 u^\lambda_\mu(z) P^C_\lambda(\vectx).
\end{equation}
The proof of Theorem~\ref{thm:Pieri} consists of two steps: 
firstly we express $u^\lambda_\mu(z)$ as a determinant;
then we interpret $u^\lambda_\mu(z)$ in terms of non-intersecting lattice paths.

The following lemma gives an explicit determinant expression for $u^\lambda_\mu(z)$.

\begin{lemma}
\label{lem:Pieri=det}
We define formal power series $b^r_s(z)$ by the relation
\begin{equation}
\label{eq:f-expansion}
\tilde{f}_s(x)
\cdot
\tilde{\Pi}_z(x)
 =
\sum_{r=0}^\infty b^r_s(z) \tilde{f}_r(x),
\end{equation}
where $\tilde{f}_d(x)$ are the Laurent polynomials defined by (\ref{eq:f}), 
and
$$
\tilde{\Pi}_z(x)
 = 
\frac{ (1 + x z) (1 + x^{-1} z) }
     { (1 - x z) (1 - x^{-1} z) }.
$$
Then we have
\begin{enumerate}
\item[(1)]
The coefficients $b^r_s(z)$ are explicitly given by
\begin{equation}
\label{eq:b}
b^r_s(z)
 =
\begin{cases}
 1 &\text{if $s=0$ and $r=0$,} \\
 2 z^r &\text{if $s=0$ and $r \ge 1$,} \\
 0 &\text{if $s \ge 1$ and $r=0$, } \\
 2 z^{s-r} (1+z^2) \dfrac{1-z^{2r}}{1-z^2} &\text{if $s \ge 1$ and $1 \le r \le s-1$,} \\
 2 (1+z^2) \dfrac{1-z^{2s}}{1-z^2} - 1 &\text{if $s \ge 1$ and $r=s$,} \\
 2 z^{r-s} (1+z^2) \dfrac{1-z^{2s}}{1-z^2} &\text{if $s \ge 1$ and $r \ge s+1$.}
\end{cases}
\end{equation}
\item[(2)]
For two sequences $\alpha = (\alpha_1, \dots, \alpha_r)$ and $\beta = (\beta_1, \dots, \beta_r)$ 
of nonnegative integers, we put $B^\alpha_\beta = \Big( b^{\alpha_i}_{\beta_j}(z) \Big)_{1 \le i, j \le r}$.
Then we have
\begin{equation}
\label{eq:u=det}
u^\lambda_\mu(z)
 =
\begin{cases}
 \det B^\lambda_\mu &\text{if $l(\lambda) = l(\mu)$,} \\
 \det B^\lambda_{\mu^0} &\text{if $l(\lambda) = l(\mu)+1$,} \\
 0 &\text{otherwise,}
\end{cases}
\end{equation}
where $\mu^0 = (\mu_1, \dots, \mu_{l(\mu)}, 0)$.
\end{enumerate}
\end{lemma}

\begin{demo}{Proof}
(1)
By a straightforward computation, we have
\begin{gather}
\label{eq:GF_f}
1 + 2 \sum_{d=0}^\infty \tilde{f}_d(x) z^r
 =
\tilde{\Pi}_z(x),
\\
\label{eq:product_f}
\tilde{f}_d(x) \tilde{f}_1(x)
 =
\begin{cases}
 \tilde{f}_1(x) &\text{if $d=0$,} \\
 \tilde{f}_2(x) &\text{if $d=1$,} \\
 \tilde{f}_{d+1}(x) + \tilde{f}_{d-1}(x) &\text{if $d \ge 2$.}
\end{cases}
\end{gather}

We proceed by induction on $s$ to prove (\ref{eq:b}).
The case $s=0$ is immediate from (\ref{eq:GF_f}).
If $s=1$, then by using (\ref{eq:GF_f}) and (\ref{eq:product_f}), we have
$$
\tilde{f}_1(x) \tilde{\Pi}_z(x) 
 =
\tilde{f}_1(x) \left( 1 + 2 \sum_{d \ge 1} \tilde{f}_d(x) z^d \right)
 =
\tilde{f}_1(x) + 2 \tilde{f}_2(x) z + \sum_{d \ge 2} (\tilde{f}_{d+1}(x) + \tilde{f}_{d-1}(x)) z^d, 
$$
from which the case $s=1$ follows.
If $s=2$, then by using $\tilde{f}_2(x) = \tilde{f}_1(x)^2$, we have
$$
\tilde{f}_2(x) \tilde{\Pi}_z(x)
 =
\tilde{f}_1(x) \left( \tilde{f}_1(x) \tilde{\Pi}_z(x) \right)
 =
b^0_1(z) \tilde{f}_1(x) + b^1_1(z) \tilde{f}_2(x)
 + \sum_{r \ge 2} b^r_1(z) \left( \tilde{f}_{r+1}(x) + \tilde{f}_{r-1}(x) \right).
$$
Hence we have
$$
b^0_2(z) = 0,
\quad
b^r_2(z) = b^{r-1}_1(z) + b^{r+1}_1(z) \quad(r \ge 1),
$$
and obtain the case $s=2$ by using the induction hypothesis.
If $s \ge 3$, then by using $\tilde{f}_s(x) = \tilde{f}_1(x) \tilde{f}_{s-1}(x) - \tilde{f}_{s-2}(x)$, we have
\begin{align*}
&
\tilde{f}_s(x) \tilde{\Pi}_z(x)
\\
 &\quad
=
\tilde{f}_1(x) \left( \tilde{f}_{s-1}(x) \tilde{\Pi}_z(x) \right) - \tilde{f}_{s-2}(x) \tilde{\Pi}_z(x)
\\
 &\quad
=
b^0_{s-1}(z) \tilde{f}_1(x) + b^1_{s-1}(z) \tilde{f}_2(x)
 + \sum_{r \ge 2} b^r_{s-1}(z) \left( \tilde{f}_{r-1}(x) + \tilde{f}_{r+1}(x) \right)
 - \sum_{r \ge 0} b^r_{s-2}(z) \tilde{f}_r(x).
\end{align*}
Hence we have
$$
b^r_s(z)
 =
\begin{cases}
-b^0_{s-2}(z) &\text{if $r=0$,} \\
b^{r-1}_{s-1}(z) + b^{r+1}_{s-1}(z) - b^r_{s-2}(z) &\text{if $r \ge 1$,}
\end{cases}
$$
and use the induction hypothesis to complete the proof of (1).

(2) can be proved in the same way as in the proof of 
\cite[Theorem~5.3 and Corollary~5.5]{Okada20}, so we omit the proof.
\end{demo}

Next we interpret the determinant in $\det B^\lambda_{\mu^*}$ with $\mu^* = \mu$ or $\mu^0$ 
as the generating function of non-intersecting lattice paths.
Let $G'$ be the directed graph with vertex set
$$
V'
 = 
\{ A_i = (i,0) : i \ge 0 \} \sqcup \{ B_i = (i,1) : i \ge 0 \} \sqcup \{ C_i = (i,2) : i \ge 0 \}
$$
and directed edge set
\begin{multline*}
E' =
\{ (A_{i+1},A_i) : i \ge 0 \} \sqcup \{ (A_i,B_i) : i \ge 0 \} \sqcup \{ (A_{i+1}, B_i) : i \ge 1 \} 
\\
\sqcup \{ (B_i, C_i) : i \ge 0 \} \sqcup \{ (B_i, C_{i+1}) : i \ge 0 \} \sqcup \{ (C_i,C_{i+1}) : i \ge 0 \}.
\end{multline*}
\begin{figure}[ht]
\centering
\setlength\unitlength{1.5pt}
\begin{picture}(110,60)(-10,-10)
\put(0,0){\circle*{2}}
\put(20,0){\circle*{2}}
\put(40,0){\circle*{2}}
\put(60,0){\circle*{2}}
\put(80,0){\circle*{2}}
\put(0,20){\circle*{2}}
\put(20,20){\circle*{2}}
\put(40,20){\circle*{2}}
\put(60,20){\circle*{2}}
\put(80,20){\circle*{2}}
\put(0,40){\circle*{2}}
\put(20,40){\circle*{2}}
\put(40,40){\circle*{2}}
\put(60,40){\circle*{2}}
\put(80,40){\circle*{2}}
\put(85,-5){\makebox(10,10){$\cdots$}}
\put(85,15){\makebox(10,10){$\cdots$}}
\put(85,35){\makebox(10,10){$\cdots$}}
\put(20,0){\vector(-1,0){19}}
\put(40,0){\vector(-1,0){19}}
\put(60,0){\vector(-1,0){19}}
\put(80,0){\vector(-1,0){19}}
\put(0,0){\vector(0,1){19}}
\put(20,0){\vector(0,1){19}}
\put(40,0){\vector(0,1){19}}
\put(60,0){\vector(0,1){19}}
\put(80,0){\vector(0,1){19}}
\put(40,0){\vector(-1,1){19}}
\put(60,0){\vector(-1,1){19}}
\put(80,0){\vector(-1,1){19}}
\put(0,20){\vector(0,1){19}}
\put(20,20){\vector(0,1){19}}
\put(40,20){\vector(0,1){19}}
\put(60,20){\vector(0,1){19}}
\put(80,20){\vector(0,1){19}}
\put(0,20){\vector(1,1){19}}
\put(20,20){\vector(1,1){19}}
\put(40,20){\vector(1,1){19}}
\put(60,20){\vector(1,1){19}}
\put(0,40){\vector(1,0){19}}
\put(20,40){\vector(1,0){19}}
\put(40,40){\vector(1,0){19}}
\put(60,40){\vector(1,0){19}}
\put(-5,40){\makebox(10,10){$C_0$}}
\put(15,40){\makebox(10,10){$C_1$}}
\put(35,40){\makebox(10,10){$C_2$}}
\put(55,40){\makebox(10,10){$C_3$}}
\put(75,40){\makebox(10,10){$C_4$}}
\put(-10,20){\makebox(10,10){$B_0$}}
\put(10,20){\makebox(10,10){$B_1$}}
\put(30,20){\makebox(10,10){$B_2$}}
\put(50,20){\makebox(10,10){$B_3$}}
\put(70,20){\makebox(10,10){$B_4$}}
\put(-5,-10){\makebox(10,10){$A_0$}}
\put(15,-10){\makebox(10,10){$A_1$}}
\put(35,-10){\makebox(10,10){$A_2$}}
\put(55,-10){\makebox(10,10){$A_3$}}
\put(75,-10){\makebox(10,10){$A_4$}}
\end{picture}
\caption{Graph $G'$}
\label{fig:G'}
\end{figure}
Note that there is no directed edge from $A_1$ to $B_0$.
Recall that a lattice path on $G'$ is a sequence $(S_0, S_1, \dots, S_m)$ of vertices of $G'$ 
such that $(S_k, S_{k+1}) \in E'$ for $0 \le k \le m-1$, 
and that a family $(P_1, \dots, P_l)$ of lattice paths on $G'$ is non-intersecting 
if no two of them have a vertex in common.
We denote by $\LL(s;r)$ the set of all lattice paths from $A_s$ to $C_r$.
For two strict partitions $\lambda$ and $\mu$ such that $l(\lambda) = l(\mu)$ or $l(\mu)+1$,
let $\LL_0(\mu;\lambda)$ be the set of all non-intersecting lattice paths 
with starting points $(A_{\mu_1}, \dots, A_{\mu_l})$ 
and ending points $(C_{\lambda_1}, \dots, C_{\lambda_l})$, where $l = l(\lambda)$.
Note that, if $l(\mu) = l-1$, then we put $\mu_l = 0$.

We define a weight function $\wt : E' \to \Int[z]$ by assigning $1$ to the vertical edges and 
$z$ to other edges:
\begin{gather*}
\wt (A_{i+1},A_i) = z,
\quad
\wt (A_i,B_i) = 1,
\quad
\wt (A_{i+1}, B_i) = z,
\\
\wt (B_i, C_i) = 1,
\quad
\wt (B_i, C_{i+1}) = z,
\quad
\wt (C_i,C_{i+1}) = z.
\end{gather*}
The weight $\wt(P)$ of a lattice path $P$ is the product of the weights of its individual steps, 
and $\wt(P_1, \dots, P_l) = \prod_{i=1}^l \wt(P_i)$.
The following lemma provides a combinatorial expression of the multiplicity $u^\lambda_\mu(z)$.

\begin{lemma}
\label{lem:Pieri=GF_path}
If $\lambda$ and $\mu$ are two strict partitions such that $l(\lambda) = l(\mu)$ or $l(\mu)+1$, 
then we have
\begin{equation}
\label{eq:u=GF_path}
u^\lambda_\mu(z)
 =
\sum_{(P_1, \dots, P_l) \in \LL_0(\mu;\lambda)} \wt(P_1, \dots, P_l),
\end{equation}
where $l = l(\lambda)$.
\end{lemma}

\begin{demo}{Proof}
We put
$$
w^s_r(z) = \sum_{P \in \Par(s;r)} \wt(P),
\quad
w^\lambda_\mu(z) = \sum_{P \in \Par_0(\mu;\lambda)} \wt(P).
$$
Then by applying the Lindstr\"om--Gessel--Vienot lemma, we have
$$
w^\lambda_\mu(z) = \det \Big( w^{\lambda_i}_{\mu_j}(z) \Big)_{1 \le i, j \le l}.
$$
Hence by Lemma~\ref{lem:Pieri=det} (2) 
it is enough to show $w^r_s(z) = b^r_s(z)$ for $r \ge 1$ and $s \ge 0$.

We proceed by induction on $r+s$.
It is easy to see that $w^1_0 = 2 z = b^1_0$ and $w^1_1 = 1+2z^2 = b^1_1$.
Suppose that $r+s \ge 2$.
Then we can show the following recurrence:
\begin{equation}
\label{eq:GF_rec}
w^r_s
 =
\begin{cases}
 z w^r_{s-1} &\text{if $r \le s-2$,} \\
 z w^{s-1}_{s-1} + z &\text{if $r=s-1$,} \\
 z^2 w^{s-1}_{s-1} + 1 + 3 z^2 &\text{if $r=s$,} \\
 z w^s_s + z &\text{if $r=s+1$,} \\
 z w^{r-1}_s &\text{if $r \ge s+2$.}
\end{cases}
\end{equation}
In fact, if $r \le s-2$, then every path in $\LL(s;r)$ passes through $A_{s-1}$, thus we have 
$w^r_s = z w^r_{s-1}$.
If $r=s-1$, then we have
$$
\LL(s;s-1)
 = 
\{ (A_s,A_{s-1})*P : P \in \LL(s-1;s-1) \} 
\cup
\{ (A_s,B_{s-1},C_{s-1}) \},
$$
where $*$ is the concatenation of paths.
Hence we have $w^{s-1}_s = z w^{s-1}_{s-1} + z$.
If $r=s$, then we have
\begin{multline*}
\LL(s;s)
 =
\{ (A_s,A_{s-1}) * P * (C_{s-1},C_s) : P \in \LL(s-1;s-1) \}
\\
\cup
\{ (A_s, A_{s-1}, B_{s-1}, C_s), (A_s, B_{s-1}, C_{s-1}, C_s), (A_s, B_{s-1}, C_s), (A_s, B_s, C_s) \},
\end{multline*}
so that $w^s_s = z^2 w^{s-1}_{s-1} + 1 + 3 z^2$.
The other cases of (\ref{eq:GF_rec}) can be checked similarly.
Then, by using the recurrence (\ref{eq:GF_rec}) and the induction hypothesis, we obtain $w^r_s = b^r_s$.
\end{demo}

Now we can complete the proof of Theorem~\ref{thm:Pieri}.

\begin{demo}{Proof of Theorem~\ref{thm:Pieri}}
By Lemma~\ref{lem:Pieri=det} (2) and (\ref{eq:GF_Pieri}), 
we have $\tilde{f}^\lambda_{\mu,(r)} = 0$ unless $l(\lambda) = l(\mu)$ or $l(\mu)+1$.

Assume from now that $l(\lambda) = l(\mu)$ or $l(\mu)+1$ and write $l = l(\lambda)$.
For a strict partition $\kappa$ of length $l$ or $l-1$, 
let $\LL_0(\mu;\lambda)_\kappa$ be the set of non-intersecting paths $(P_1, \dots, P_l) 
\in \LL_0(\mu;\lambda)$ such that $P_i$ passes through the vertex $B_{\kappa_i}$ for $1 \le i \le l$.
Then we have
$$
\LL_0(\mu;\lambda) = \bigsqcup_\kappa \LL_0(\mu;\lambda)_\kappa,
$$
where $\kappa$ runs over all strict partitions of length $l$ or $l-1$.

First we show that $\LL_0(\mu;\lambda)_\kappa = \emptyset$ unless $\mu \succ \kappa$ and $\lambda \succ \kappa$.
Suppose that there exists a family of non-intersecting lattice paths 
$(P_1, \dots, P_l) \in \LL_0(\mu;\lambda)_\kappa$.
Since the $i$th path $P_i$ starts at $A_{\mu_i}$ and passes through $B_{\kappa_i}$, 
we see that $\mu_i \ge \kappa_i$.
If $\mu_i > \kappa_i$, then $P_i$ passes through $A_{\kappa_i+1}$ and does not 
intersect with $P_{i+1}$, hence $\kappa_i+1 > \mu_{i+1}$, i.e., $\kappa_i \ge \mu_{i+1}$.
Similarly we can check $\lambda \succ \kappa$.

Next it is clear from the definition of the weight that, if $(P_1, \dots, P_l) \in \LL_0(\mu;\lambda)_\kappa$, then we have
$$
\wt (P_i) = z^{(\mu_i-\kappa_i)+(\lambda_i-\kappa_i)},
\quad
\wt (P_1, \dots, P_l)
 =
z^{(|\mu|-|\kappa|) + (|\lambda|-|\kappa|)}.
$$
Hence, by Lemma~\ref{lem:Pieri=GF_path} and (\ref{eq:GF_Pieri}), we have
$$
2 \tilde{f}^\lambda_{\mu,(r)}
 =
\sum_\kappa \# \LL_0(\mu;\lambda)_\kappa,
$$
where $r>0$ and $\kappa$ runs over all strict partitions satisfying $\mu \succ \kappa$, $\lambda \succ \kappa$ 
and $(|\mu|-|\kappa|) + (|\lambda|-|\kappa|) = r$.

Finally we compute the cardinality of $\LL_0(\mu;\lambda)_\kappa$.
Suppose that $\mu \succ \kappa$ and $\lambda \succ \kappa$.
If $1 \le i \le l-1$, we have
\begin{gather*}
\# \{ P \in \LL(\mu_i;\kappa_i) : \text{$P$ does not pass through $A_{\mu_{i+1}}$} \}
 =
\begin{cases}
 1 &\text{if $\kappa_i = \mu_i$,} \\
 2 &\text{if $\mu_{i+1} < \kappa_i < \mu_i$,} \\
 1 &\text{if $\kappa_i = \mu_{i+1}$,}
\end{cases}
\\
\# \{ P \in \LL(\kappa_i;\lambda_i) : \text{$P$ does not pass through $A_{\lambda_{i+1}}$} \}
 =
\begin{cases}
 1 &\text{if $\kappa_i = \lambda_i$,} \\
 2 &\text{if $\lambda_{i+1} < \kappa_i < \lambda_i$,} \\
 1 &\text{if $\kappa_i = \lambda_{i+1}$,}
\end{cases}
\end{gather*}
Also we have
$$
\# \LL(\mu_l;\kappa_l)
 =
\begin{cases}
 1 &\text{if $\mu_l > 0$ and $\kappa_l = \mu_l$,} \\
 2 &\text{if $\mu_l > 0$ and $0 < \kappa_l < \mu_l$,} \\
 1 &\text{if $\mu_l > 0$ and $\kappa_l = 0$,} \\
 1 &\text{if $\mu_l = 0$ and $\kappa_l = 0$,}
\end{cases}
\quad
\# \LL(\lambda_l;\kappa_l)
 =
\begin{cases}
 1 &\text{if $\kappa_l = \lambda_l$,} \\
 2 &\text{if $\kappa_l < \lambda_l$.}
\end{cases}
$$
Note that $\lambda_l > 0$.
Hence it follows from (\ref{eq:a}) that
$$
\# \LL_0(\mu;\lambda)_\kappa
 =
2^{a(\lambda,\kappa) + a(\mu,\kappa) - \chi[l(\mu)>l(\kappa)]}.
$$
This completes the proof of Theorem~\ref{thm:Pieri}.
\end{demo}
\section{%
Positivity conjectures
}

In this section, we present some positivity conjectures on various expansion coefficients 
involving symplectic $P$-functions. 
The conjectures in this section have been checked on a computer 
with a help of the Maple package {\tt SF} developed by Stembridge \cite{SF}.
It would be interesting to resolve these conjecture 
by finding combinatorial rules for describing the coefficients.

\subsection{%
Structure constants for multiplication
}

Let $\Lambda$ be the ring of symmetric functions and $\Gamma$ the subring generated by $q_r$ ($r \ge 1$).
Since the Schur $P$-functions $P_\lambda$ form a basis of $\Gamma$, 
we can expand the product $P_\mu \cdot P_\nu$ in the form
\begin{equation}
\label{eq:P-product}
P_\mu \cdot P_\nu
 =
\sum_{\lambda \in \SPar} f^\lambda_{\mu,\nu} P_\lambda,
\end{equation}
where $\SPar$ is the set of strict partitions.
Since Schur $P$-functions are homogeneous, 
we have $f^\lambda_{\mu,\nu} = 0$ unless $|\lambda| = |\mu| + |\nu|$.
It is known that the coefficients $f^\lambda_{\mu,\nu}$ are nonnegative integers, 
and there are several combinatorial models for the coefficients $f^\lambda_{\mu,\nu}$ 
such as \cite[Sectioin~7.2]{Worley84}, \cite[Theorem~8.3]{Stembridge89}, \cite[Corollary~5.14]{Cho13} 
and \cite[Theorem~4.13]{GJKKK14}.

By Proposition~\ref{prop:uQC_basis} (2), 
the universal symplectic $P$-functions $\Pfun^C_\lambda$ 
form another basis of $\Gamma$. 
Given three strict partitions $\lambda$, $\mu$ and $\nu$, 
the structure constant $\tilde{f}^\lambda_{\mu,\nu}$ is defined as the coefficient of $\Pfun^C_\lambda$ 
in the expansion
\begin{equation}
\label{eq:uPC-product}
\Pfun^C_\mu \cdot \Pfun^C_\nu
 =
\sum_{\lambda \in \SPar} \tilde{f}^\lambda_{\mu,\nu} \Pfun^C_\lambda.
\end{equation}
Since $\Gamma$ is a commutative ring with unit $1 = \Pfun^C_\emptyset$, 
we have $\tilde{f}^\lambda_{\mu,\nu} = \tilde{f}^\lambda_{\nu,\mu}$ 
and $\tilde{f}^\lambda_{\mu,\emptyset} = \tilde{f}^\lambda_{\emptyset,\mu} = \delta_{\mu,\nu}$.
And it follows from (\ref{eq:uQC_by_Q}) 
that, if $|\lambda| = |\mu|+|\nu|$, then $\tilde{f}^\lambda_{\mu,\nu} = f^\lambda_{\mu,\nu} \ge 0$, 
and that $\tilde{f}^\lambda_{\mu,\nu} = 0$ unless $|\lambda| \le |\mu|+|\nu|$ and $|\mu|+|\nu|-|\lambda|$ is even.

By applying the ring homomorphism $\tilde{\pi}_n : \Lambda \to \Rat[x_1^{\pm 1}, \dots, x_n^{\pm 1}]^{W_n}$ 
given by (\ref{eq:pi}) to the both sides of (\ref{eq:uPC-product}), we have
\begin{equation}
\label{eq:PC-product}
P^C_\mu (x_1, \dots, x_n) \cdot P^C_\nu (x_1, \dots, x_n)
 =
\sum_\lambda \tilde{f}^\lambda_{\mu,\nu} P^C_\lambda(x_1, \dots, x_n),
\end{equation}
where $\mu$ and $\nu$ are strict partitions of length $\le n$ and 
$\lambda$ runs over all strict partitions of length $\le n$.

\begin{conjecture}
\label{conj:1}
The structure constants $\tilde{f}^\lambda_{\mu,\nu}$ for the multiplication of symplectic $P$-functions 
are nonnegative integers.
\end{conjecture}

We have checked on a computer the validity of this conjecture for all pairs $(\mu, \nu)$ of strict partitions 
satisfying $|\mu| + |\nu| \le 20$.
And, in support of this conjecture, we have the following evidence.

\begin{prop}
\label{prop:conj1}
Conjecture~\ref{conj:1} holds in the following cases:
\begin{enumerate}
\item[(1)]
$|\lambda| = |\mu| + |\nu|$.
In this case, $\tilde{f}^\lambda_{\mu,\nu} = f^\lambda_{\mu,\nu} \ge 0$.
\item[(2)]
$l(\mu)=1$ or $l(\nu)=1$.
In this case, the coefficients $\tilde{f}^\lambda_{\mu,(r)}$ are given by the Pieri rule 
in Theorem~\ref{thm:Pieri}.
\item[(3)]
$\mu = \delta_r = (r,r-1,\dots,2,1)$ and $\nu = \delta_s = (s,s-1,\dots,2,1)$.
In this case, we have
$$
\Pfun^C_{\delta_r} \cdot \Pfun^C_{\delta_s} = \Pfun^C_{\delta_r+\delta_s}.
$$
\end{enumerate}
\end{prop}

\begin{demo}{Proof}
It remains to prove (3).
By performing row/column operations on the Pfaffian in the numerator of the Nimmo-type formula (\ref{eq:Nimmo_PC}) 
with use of the relation 
$$
(x+x^{-1})^d
 =
\sum_{i=0}^{\lfloor d/2 \rfloor}
 \left\{ \binom{d-2}{i} - \binom{d-2}{i-2} \right\} \tilde{g}_{d-2i}(x),
$$
we can show
\begin{equation}
\label{eq:PC_delta}
P^C_{\delta_r+\delta_s}(\vectx)
 =
P_{\delta_r+\delta_s}(\vectx+\vectx^{-1}),
\end{equation}
where $\vectx = (x_1, \dots, x_n)$ and $\vectx+\vectx^{-1} = (x_1+x_1^{-1}, \dots, x_n+x_n^{-1})$.
In particular, $P^C_{\delta_r}(\vectx) = P_{\delta_r}(\vectx+\vectx^{-1})$.

Since $P_{\delta_r} P_{\delta_s} = P_{\delta_r+\delta_s}$ 
(see \cite[(7.17)]{Worley84} for a combinatorial proof 
and \cite[Theorem~5.2]{Okada14} for a sketch of an algebraic proof), 
we obtain
$$
P^C_{\delta_r}(\vectx) P^C_{\delta_s}(\vectx)
 =
P_{\delta_r}(\vectx+\vectx^{-1}) P_{\delta_s}(\vectx+\vectx^{-1})
 =
P_{\delta_r+\delta_s}(\vectx+\vectx^{-1})
 =
P^C_{\delta_r+\delta_s}(\vectx).
$$
Hence by Lemma~\ref{lem:proj} (1) 
we conclude $\Pfun^C_{\delta_r} \cdot \Pfun^C_{\delta_s} = \Pfun^C_{\delta_r+\delta_s}$ in $\Lambda$.
\end{demo}

\subsection{%
Structure constants for comultiplication
}

The ring of symmetric functions $\Lambda$ has a structure of graded Hopf algebra, 
whose coproduct $\Delta$ is given by the ``alphabet doubling trick''.
For $f \in \Lambda$, we have $\Delta(f) = \sum_{i=1}^k g_i \otimes h_i$ if and only if 
$f(X \cup Y) = \sum_{i=1}^k g_i(X) h_i(Y)$, where $X$ and $Y$ are two disjoint sets of infinitely many variables.
Since $q_r(X \cup Y) = \sum_{k=0}^r q_k(X) q_{r-k}(Y)$ by (\ref{eq:GF_q}), 
we see that $\Gamma$ is the Hopf subalgebra of $\Lambda$.
For Schur $Q$-functions, we have
\begin{equation}
\label{eq:Q-coproduct}
Q_\lambda (X \cup Y)
 =
\sum_{\mu,\nu \in \SPar} f^\lambda_{\mu,\nu} Q_\mu(X) Q_\nu(Y),
\end{equation}
where $f^\lambda_{\mu,\nu}$ are the structure constants for the multiplication of Schur $P$-functions 
given by (\ref{eq:P-product}).
Note that (\ref{eq:Q-coproduct}) is equivalent to $Q_{\lambda/\mu} = \sum_\nu f^\lambda_{\mu,\nu} Q_\nu$ 
for skew $Q$-functions.

Let $\tilde{d}^\lambda_{\mu,\nu}$ be the structure constant for comultiplication of symplectic $Q$-functions 
defined by
\begin{equation}
\label{eq:uQC-coproduct}
\Qfun^C_\lambda(X \cup Y)
 =
\sum_{\mu,\nu \in \SPar} \tilde{d}^\lambda_{\mu,\nu} \Qfun^C_\mu(X) \Qfun^C_\nu(Y).
\end{equation}
By Proposition~\ref{prop:sep_var}, we see that (\ref{eq:uQC-coproduct}) is equivalent to
$$
\Qfun^C_{\lambda/\mu} = \sum_{\nu \in \SPar} \tilde{d}^\lambda_{\mu,\nu} \Qfun^C_\nu.
$$
It follows from (\ref{eq:uQC_by_Q}) 
that, if $|\lambda| = |\mu|+|\nu|$, then $\tilde{d}^\lambda_{\mu,\nu} = f^\lambda_{\mu,\nu} \ge 0$, 
and that $\tilde{d}^\lambda_{\mu,\nu} = 0$ unless $|\mu|+|\nu| \le |\lambda|$ and $|\lambda| - |\mu|-|\nu|$ is even.
Since $\Qfun^C_{\lambda,\emptyset} = \Qfun^C_\lambda$, 
we have $\tilde{d}^\lambda_{\mu,\emptyset} = \tilde{d}^\lambda_{\emptyset,\mu} = \delta_{\lambda,\mu}$.
And, since $\Qfun^C_{\lambda/\mu} = 0$ unless $\lambda \supset \mu$ (Proposition~\ref{prop:uskewQC}), 
we see that $\tilde{d}^\lambda_{\mu,\nu} = 0$ unless $\lambda \supset \mu$ and $\lambda \supset \nu$.
By applying $\tilde{\pi}_n \otimes \tilde{\pi}_m$ to the both sides of (\ref{eq:uQC-coproduct}), we obtain
$$
Q^C_\lambda(x_1, \dots, x_n, y_1, \dots, y_m)
 =
\sum_{\mu,\nu} \tilde{d}^\lambda_{\mu,\nu} Q^C_\mu(x_1, \dots, x_n) Q^C_\nu(y_1, \dots, y_m),
$$
where $\lambda$ is a strict partition of length $\le n+m$ 
and the sum is taken over all pairs $(\mu,\nu)$ of strict partitions such that $l(\mu) \le n$ and $l(\nu) \le m$.

\begin{conjecture}
\label{conj:2}
The structure constants $\tilde{d}^\lambda_{\mu,\nu}$ for the comultiplication of symplectic $Q$-functions 
are nonnegative integers.
\end{conjecture}

We have checked that this conjecture holds for strict partitions $\lambda$ with $|\lambda| \le 20$.
And we can prove the following proposition.

\begin{prop}
\label{prop:conj2}
Conjecture~\ref{conj:2} holds in the following cases.
\begin{enumerate}
\item[(1)]
$|\lambda| = |\mu| + |\nu|$.
In this case, we have $\tilde{d}^\lambda_{\mu,\nu} = f^\lambda_{\mu,\nu} \ge 0$.
\item[(2)]
$l(\lambda) = 1$.
In this case, we have 
$\Qfun^C_{(r)}(X \cup Y) = \sum_{k=0}^r \Qfun^C_{(k)}(X) \Qfun^C_{(r-k)}(Y)$.
\item[(3)]
$\lambda = \delta_r$.
In this case, we have 
$\Qfun^C_{\delta_r}(X \cup Y) = \sum_\mu \Qfun^C_\mu(X) \Qfun^C_{\mu^*}(Y)$, 
where $\mu$ runs over all strict partitions such that $\mu \subset \delta_r$ 
and $\mu^*$ is the complement of $\mu$ with respect to $\delta_r$.
\end{enumerate}
\end{prop}

\begin{demo}{Proof}
(1) is already discussed above.
(2) is a consequence of (\ref{eq:GF_q}) since $\Qfun^C_{(r)} = q_r$ by (\ref{eq:uQC_length1}).
(3) follows from $\Qfun^C_{\delta_r/\mu} = \Qfun^C_{\mu^*}$ (see Proposition~\ref{prop:flip}).
\end{demo}

\subsection{
Symplectic $Q$-functions and symplectic Schur functions
}

Since Schur $P$-functions are symmetric functions, we can expand them in the Schur function basis:
\begin{equation}
\label{eq:P_by_s}
P_\lambda
 =
\sum_{\mu \in \Par} g_{\lambda,\mu} s_\mu,
\end{equation}
where $\Par$ is the set of partitions.
Then it is known that the expansion coefficients $g_{\lambda,\mu}$ are nonnegative integers, 
and there are several combinatorial models for the coefficients $g_{\lambda,\mu}$ 
such as \cite[(7.12)]{Worley84}, \cite[Theorem~13.1]{Sagan87} and \cite[Theorem~9.3]{Stembridge89}.
In this subsection we consider the expansion of universal symplectic $P$-functions 
into universal symplectic Schur functions.

The \emph{universal symplectic Schur function} $s^C_\lambda \in \Lambda$ 
corresponding to a partition $\lambda$ is defined by
\begin{equation}
\label{eq:uSC}
s^C_\lambda
 =
\frac{ 1 }{ 2 }
\det \Big(
 h_{\lambda_i-i+j} + h_{\lambda_i-i-j+2}
\Big)_{1 \le i, j \le l(\lambda)},
\end{equation}
where $h_d$ is the $d$th complete symmetric function.
For a partition $\lambda$ of length $\le n$, 
we have $\tilde{\pi}_n( s^C_\lambda ) = S^C_\lambda(x_1, \dots, x_n)$, 
so $s^C_\lambda$ is a symmetric function lift of the irreducible characters of symplectic groups 
corresponding to $\lambda$.
See \cite{KT87}, \cite{Koike89} and \cite{King90} for expositions 
of classical group characters and their lifts to symmetric functions, 
called universal characters.

It follows from (\ref{eq:uSC}) that $s^C_\lambda$ can be written 
as a linear combination of Schur functions in the form
\begin{equation}
\label{eq:uSC_by_s}
s^C_\lambda = s_\lambda + \sum_\mu a'_{\lambda,\mu} s_\mu,
\end{equation}
where $\mu$ runs over all partitions such that $|\mu| < |\lambda|$ and $|\lambda|-|\mu|$ is even.
Hence the universal symplectic Schur functions $\{ s^C_\lambda : \lambda \in \Par \}$ 
form a basis of $\Lambda$. 
We can define the expansion coefficients $\tilde{g}_{\lambda,\mu}$ by the relation
\begin{equation}
\label{eq:uPC_by_uSC}
\Pfun^C_\lambda
 =
\sum_{\mu \in \Par} \tilde{g}_{\lambda,\mu} s^C_\mu,
\end{equation}
where $\lambda$ is a strict partition.
By using (\ref{eq:uQC_by_Q}) and (\ref{eq:uSC_by_s}), 
we can see that 
if $|\lambda| = |\mu|$, then $\tilde{g}_{\lambda,\mu} = g_{\lambda,\mu} \ge 0$, 
and that $\tilde{g}_{\lambda,\mu} = 0$ unless $|\lambda| \ge |\mu|$ and $|\lambda|-|\mu|$ is even.

\begin{conjecture}
\label{conj:3-1}
The expansion coefficients $\tilde{g}_{\lambda,\mu}$ in (\ref{eq:uPC_by_uSC}) are nonnegative integers.
\end{conjecture}

Since the symplectic Schur functions $\{ S^C_\lambda(x_1, \dots, x_n) : \lambda \in \Par^{(n)} \}$ 
form a basis of $\Rat[x_1^{\pm 1}, \dots, x_n^{\pm 1}]^{W_n}$, 
where $\Par^{(n)}$ is the set of partitions of length $\le n$, 
we can define $\tilde{g}_{\lambda,\mu}(n)$ by the relation
\begin{equation}
\label{eq:PC_by_SC}
P^C_\lambda(x_1, \dots, x_n)
 =
\sum_{\mu \in \Par^{(n)}} \tilde{g}_{\lambda,\mu}(n) S^C_\mu(x_1, \dots, x_n),
\end{equation}
where $\lambda$ is a strict partition of length $\le n$.

\begin{conjecture}
\label{conj:3-2}
For a positive integer $n$, 
the expansion coefficients $\tilde{g}_{\lambda,\mu}(n)$ in (\ref{eq:PC_by_SC}) are nonnegative integers.
\end{conjecture}

Since the specialization $\tilde{\pi}_n(s^C_\mu)$ is not a nonnegative linear combination 
of symplectic Schur functions in general, 
Conjecture~\ref{conj:3-1} does not imply Conjecture~\ref{conj:3-2}.
On the other hand, if $n$ is large enough for a given partition $\lambda$ (e.g. $n \ge |\lambda|$), 
then $\tilde{\pi}_n(s^C_\mu) = S^C_\mu(x_1, \dots, x_n)$ for any partitions $\mu$ such that $|\mu| \le |\lambda|$, 
and thus we have $\tilde{g}_{\lambda,\mu}(n) = \tilde{g}_{\lambda,\mu}$.

By multiplying the both sides of (\ref{eq:PC_by_SC}) with $P^C_{\delta_n}(x_1, \dots, x_n)$ 
and using (\ref{eq:PC=SC*SC}), we have
$$
P^C_\lambda(x_1, \dots, x_n) \cdot P^C_{\delta_n}(x_1, \dots, x_n)
 =
\sum_{\mu \in \Par^{(n)}} \tilde{g}_{\lambda,\mu}(n) P^C_{\mu+\delta_n}(x_1, \dots, x_n).
$$
By comparing this with (\ref{eq:PC-product}), we see that 
$\tilde{g}_{\lambda,\mu}(n) = \tilde{f}^{\mu+\delta_n}_{\lambda,\delta_n}$ 
if $l(\lambda)$, $l(\mu) \le n$.
Hence Conjecture~\ref{conj:3-2} is a consequence of Conjecture~\ref{conj:1}.

We have verified Conjecture~\ref{conj:3-1} for all strict partitions $\lambda$ with $|\lambda| \le 18$.
And we can prove the following special cases of Conjectures~\ref{conj:3-1} and \ref{conj:3-2}.

\begin{prop}
\label{prop:conj3-1}
Conjecture~\ref{conj:3-1} is true in the following cases:
\begin{enumerate}
\item[(1)]
$|\lambda| = |\mu|$.
In this case, we have $\tilde{g}_{\lambda,\mu} = g_{\lambda,\mu} \ge 0$.
\item[(2)]
$l(\lambda) = 1$.
In this case we have
$$
\tilde{g}_{(r),\mu}
 =
\begin{cases}
 1 &\text{if $\mu$ is a hook and $|\mu| = r$,} \\
 2 &\text{if $\mu$ is a hook, $|\mu| \equiv r \bmod 2$ and $1 \le |\mu| \le r-1$,} \\
 1 &\text{if $\mu = \emptyset$ and $r$ is even,} \\
 0 &\text{otherwise,}
\end{cases}
$$
where a partition $\mu$ is a hook if $\lambda = (a+1,1^b)$ for some nonnegative integers $a$ and $b$.
\item[(3)]
$\lambda = \delta_r + \delta_s$.
\end{enumerate}
\end{prop}

\begin{demo}{Proof}
We prove (2) and (3).

(2)
Since $\Pfun^C_{(r)} = q_r/2$ by definition, it follows from \cite[III.8 Example~6(c)]{Macdonald95} that
$$
\Pfun^C_{(r)} = \sum_{a+b = r-1} s_{(a|b)},
$$
where $(a|b) = (a+1, 1^b)$.
Here we use the following Littlewood formula (\ref{eq:s_by_uSC}) (see e.g. \cite[Theorem~2.3.1 (1)]{KT87}) 
to express $s_{(a,b)}$ in terms of universal symplectic Schur functions.
For an arbitrary partition $\lambda$, the corresponding Schur function $s_\lambda$ is expressed in the form
\begin{equation}
\label{eq:s_by_uSC}
s_\lambda
 = 
\sum_{\mu \in \Par} \left( \sum_\nu c^\lambda_{\mu,\nu} \right) s^C_\mu,
\end{equation}
where $\nu$ runs over all partitions such that all columns have even length, 
and $c^\alpha_{\beta,\gamma}$ stands for the ordinary Littlewood--Richardson coefficients.
If $\lambda$ is a hook and $\nu$ is a partition such that all columns have even length, 
then $c^\lambda_{\mu,\nu} = 0$ unless $\nu$ is a one-column partition $(1^m)$.
Hence we see that
\begin{align*}
s_{(a|b)}
 &=
\sum_{i=0}^{\lfloor b/2 \rfloor} s^C_{(a|b-2i)}
 +
\sum_{i=0}^{\lfloor (b-1)/2 \rfloor} s^C_{(a-1|b-2i-1)}
\quad(a \ge 1),
\\
s_{(0|b)}
 &=
\sum_{i=0}^{\lfloor b/2 \rfloor} s^C_{(0|b-2i)}
 +
\begin{cases}
 0 &\text{if $b$ is even,} \\
 s^C_\emptyset &\text{if $b$ is odd.}
\end{cases}
\end{align*}
By combining these relations we obtain the desired result.

(3)
The bialternant formula (\ref{eq:SC}) can be transformed into
$$
S^C_\lambda(\vectx)
 =
\frac{ 1 }
     { \prod_{1 \le i < j \le n} ( (x_i+x_i^{-1}) - (x_j+x_j^{-1}) ) }
\det \Big(
 c_{\lambda_j+n-j}(x_i)
\Big)_{1 \le i, j \le n},
$$
where $\vectx = (x_1, \dots, x_n)$ and $c_d(x) = (x^{d+1}-x^{-(d+1)})/(x-x^{-1})$.
By performing column operations on the determinant above 
with use of the relation
$$
(x+x^{-1})^d
 =
\sum_{i=0}^{\lfloor d/2 \rfloor}
 \left\{ \binom{d-1}{i} - \binom{d-1}{i-2} \right\} \tilde{c}_{d-2i}(x),
$$
we can obtain
$$
S^C_{\delta_r}(\vectx)
 =
s_{\delta_r}(\vectx+\vectx^{-1}).
$$
Hence, by using (\ref{eq:PC_delta}) and 
$P_{\delta_r+\delta_s}(\vectx) = s_{\delta_r}(\vectx) s_{\delta_s}(\vectx)$ 
(see \cite[(7.17)]{Worley84}, \cite[Theorem~5.2]{Okada14}), 
we have
$$
P^C_{\delta_r+\delta_s}(\vectx)
 =
P_{\delta_r+\delta_s}(\vectx+\vectx^{-1})
 =
s_{\delta_r}(\vectx+\vectx^{-1})
s_{\delta_s}(\vectx+\vectx^{-1})
 =
S^C_{\delta_r}(\vectx)
S^C_{\delta_s}(\vectx).
$$
Thus, by Lemma~\ref{lem:proj} (1), we see that 
$\Pfun^C_{\delta_r+\delta_s} = s^C_{\delta_r} s^C_{\delta_s}$ in $\Lambda$.
The Newell--Littlewood formula (see e.g. \cite[Theorem~3.1]{Koike89}) asserts 
that the product $s^C_\mu s^C_\nu$ is positively expanded in the universal symplectic Schur function basis.
Hence we see that $\Pfun^C_{\delta_r+\delta_s}$ 
is a nonnegative linear combination of universal symplectic Schur functions.
\end{demo}

\begin{prop}
\label{prop:conj3-2}
Conjecture~\ref{conj:3-2} is true in the following cases:
\begin{enumerate}
\item[(1)]
$l(\lambda) = 1$.
\item[(2)]
$l(\lambda) = n$.
\item[(3)]
$\lambda = \delta_r + \delta_s$ with $r$, $s \le n$.
\end{enumerate}
\end{prop}

\begin{demo}{Proof}
(1)
Let $e_r(\vectx,\vectx^{-1})$ and $h_r(\vectx,\vectx^{-1})$ be 
the $r$th elementary and complete symmetric polynomials in 
$(\vectx,\vectx^{-1}) = (x_1, \dots, x_n, x_1^{-1}, \dots, x_n^{-1})$.
Then it follows from (\ref{eq:GF_length1}) that
$$
Q^C_{(r)}(\vectx) = \sum_{p+q=r} e_p(\vectx,\vectx^{-1}) h_q(\vectx,\vectx^{-1}).
$$
Note that $e_p(\vectx,\vectx^{-1})$ and $h_q(\vectx,\vectx^{-1})$ are the characters 
of the exterior power $\twedge^p(V)$ and the symmetric powers $S^q(V)$ 
of the vector representation $V = \Comp^{2n}$ of $\Symp_{2n}(\Comp)$ respectively.
Hence $Q^C_{(r)}(\vectx)$ is the character of some representation of $\Symp_{2n}(\Comp)$, 
and a nonnegative linear combination of symplectic Schur functions.

(2)
If $l(\lambda) = n$, then $\lambda = \mu + \delta_n$ for some partition $\mu$ of length $\le n$ 
and $P^C_\lambda(\vectx) = S^C_{\delta_n}(\vectx) S^C_\mu(\vectx)$ by (\ref{eq:PC=SC*SC}).
Hence it is a nonnegative linear combination of symplectic Schur functions.

(3)
The claim follows from $P^C_{\delta_r+\delta_s}(\vectx) = S^C_{\delta_r}(\vectx) S^C_{\delta_s}(\vectx)$, 
which was proved in the proof of Proposition~\ref{prop:conj3-1}.
\end{demo}

\subsection{%
Schur $P$-functions and symplectic $P$-functions
}

It follows from (\ref{eq:uQC_by_Q}) that the Schur $P$-function $P_\lambda$ 
can be expressed as a linear combination of universal symplectic $P$-functions 
in the form
\begin{equation}
\label{eq:P_by_uPC}
P_\lambda
 = 
\Pfun^C_\lambda
 +
\sum_\mu b_{\lambda,\mu} \Pfun^C_\mu,
\end{equation}
where $\mu$ runs over all strict partitions such that $|\mu| < |\lambda|$ and $|\lambda| - |\mu|$ is even.
This expansion is a $P$-function analogue of (\ref{eq:s_by_uSC}).
By applying $\tilde{\pi}_n$, we have
\begin{equation}
\label{eq:P_by_PC}
P_\lambda(x_1, \dots, x_n, x_1^{-1}, \dots, x_n^{-1})
 =
P^C_\lambda(x_1, \dots, x_n)
 +
\sum_\mu b_{\lambda,\mu} P^C_\mu(x_1, \dots, x_n),
\end{equation}
where $\mu$ runs over all strict partitions satisfying 
$|\mu| < |\lambda|$, $|\lambda| - |\mu|$ is even and $l(\mu) \le n$.

\begin{conjecture}
\label{conj:4}
The expansion coefficients $b_{\lambda,\mu}$ in the expansion (\ref{eq:P_by_uPC}) are nonnegative integers.
\end{conjecture}

This conjecture has been checked by computer for strict partitions $\lambda$ with $|\lambda| \le 20$.
And we can prove the case where $l(\lambda) \le 2$.

\begin{prop}
\label{prop:conj4}
Conjecture~\ref{conj:4} holds in the following cases:
\begin{enumerate}
\item[(1)]
$l(\lambda) = 1$.
In this case, we have $P_{(r)} = \Pfun^C_{(r)} = q_r/2$.
\item[(2)]
$l(\lambda) = 2$.
In this case, we have
\begin{equation}
\label{eq:conj4_length2}
P_{(r,s)}
 =
\Pfun^C_{(r,s)}
 +
2 \sum_{j=1}^{s-1} \Pfun^C_{(r-j,s-j)}
 +
\Pfun^C_{(r-s)}.
\end{equation}
\end{enumerate}
\end{prop}

\begin{demo}{Proof}
(1) is obvious from the definition.

(2)
By using the definitions (\ref{eq:Q_length2}) and (\ref{eq:uQC_length2}), we have
\begin{gather*}
Q_{(r,1)} = q_r q_1 - 2 q_{r+1},
\quad
\Qfun^C_{(r,1)} = q_r q_1 - 2 q_{r+1} - 2 q_{r-1},
\\
Q_{(r,s)} - Q_{(r-1,s-1)} = \Qfun^C_{(r,s)} + \Qfun^C_{(r-1,s-1)} \quad(s \ge 2).
\end{gather*}
Now the induction on $s$ shows $Q_{(r,s)} = \Qfun^C_{(r,s)} + 2 \sum_{j=1}^s \Qfun^C_{(r-j,s-j)}$, 
which is equivalent to (\ref{eq:conj4_length2}).
\end{demo}
\section{%
Factorial symplectic $Q$-functions
}

In this section, we prove a factorial version of Theorem~\ref{thm:tableaux} 
(a tableau description of symplectic $Q$-functions), 
which implies that the factorial symplectic $Q$-functions introduced by 
Foley--King \cite{FK18} are obtained by specializing $t=-1$ in the factorial Hall--Littlewood 
functions associated to the root system of type $C_n$.

\subsection{%
Factorial symplectic $Q$-functions
}

For factorial parameters $\vecta = (a_0, a_1, a_2, \dots)$, 
the factorial monomials $(x|\vecta)^r$ ($r \ge 0$) are given by
$$
(x|\vecta)^r = \prod_{i=0}^{r-1} (x+a_i).
$$
We introduce the \emph{factorial symplectic Hall--Littlewood function} $\PP^C_\lambda(\vectx|\vecta;t)$ 
corresponding to a partition $\lambda$ of length $\le n$ 
by replacing the monomial $x_i^{\lambda_i}$ with the factorial monomial $(x_i|\vecta)^{\lambda_i}$ 
in the definition (\ref{eq:HLC}):
\begin{equation}
\label{eq:fHLC}
\PP^C_\lambda (\vectx|\vecta;t)
 =
\frac{1}{v^{(n)}_\lambda(t)}
\sum_{w \in W_n}
 w \left(
  \prod_{i=1}^n (x_i|\vecta)^{\lambda_i}
  \prod_{i=1}^n \frac{ 1 - t x_i^{-2} }{ 1 - x_i^{-2} }
  \prod_{1 \le i < j \le n}
   \frac{ 1 - t x_i^{-1} x_j }{ 1 - x_i^{-1} x_j }
   \frac{ 1 - t x_i^{-1} x_j^{-1} }{ 1 - x_i^{-1} x_j^{-1} }
 \right),
\end{equation}
where $v^{(n)}_\lambda(t)$ is given by (\ref{eq:v}).
Then, for a strict partition of length $\le n$, 
the \emph{factorial symplectic $P$-function} and \emph{factorial symplectic $Q$-function} 
are defined by
\begin{equation}
\label{eq:fPC}
P^C_\lambda(\vectx|\vecta) = \PP^C_\lambda(\vectx|\vecta;-1),
\quad
Q^C_\lambda(\vectx|\vecta) = 2^{l(\lambda)} \PP^C_\lambda(\vectx|\vecta;-1)
\end{equation}
respectively.
When the factorial parameters $a_i$ are all zero, then these functions reduce to 
the symplectic $P$-/$Q$-functions.

The aim of this section is to prove a tableau description of $Q^C_\lambda(\vectx|\vecta)$, 
which enables us to identify our factorial symplectic $Q$-functions with
Foley--King's factorial symplectic $Q$-function defined in \cite[Definition~6]{FK18}.
Following \cite{FK18}, we introduce the factorial weight $(\vectx|\vecta)^T$ 
of a symplectic primed shifted tableau $T$ of shape $\lambda/\mu$ 
by putting
$$
(\vectx|\vecta)^T
 =
\prod_{(i,j) \in S(\lambda/\mu)} w(T_{i,j};i,j),
$$
where $w(\gamma;i,j)$ is given by
$$
\wt(\gamma;i,j)
 =
\begin{cases}
x_k - a_{j-i} &\text{if $\gamma = k'$,} \\
x_k + a_{j-i} &\text{if $\gamma = k$,} \\
x_k^{-1} - a_{j-i} &\text{if $\gamma = \overline{k}'$,} \\
x_k^{-1} + a_{j-i} &\text{if $\gamma = \overline{k}$.} \\
\end{cases}
$$
The main result of this section is the following theorem.

\begin{theorem}
(See \cite[Theorem~15]{FK18})
\label{thm:fac_tableaux}
Suppose that $a_0 = 0$.
For a strict partition $\lambda$ of length $\le n$, we have
\begin{equation}
\label{eq:fac_tableaux}
Q^C_\lambda(\vectx|\vecta)
 =
\sum_{T \in \QTab^C_n(\lambda)} (\vectx|\vecta)^T.
\end{equation}
\end{theorem}

The proof is postponed to Section~7.3, 
where we extend this theorem to a skew version (Theorem~\ref{thm:fac_skew_tableaux}).
Toward the proof, we establish some formulas for $Q^C_\lambda(\vectx|\vecta)$.
The following proposition is proved by the same argument as in the proof 
of \cite[Theorem~7.2]{Okada20}, 
 
\begin{prop}
\label{prop:fQC=genP}
We define Laurent polynomials $\tilde{g}_d(x|\vecta)$ and polynomials $g_d(u|\vecta)$ by
$$
\tilde{g}_d(x|\vecta)
 =
g_d(x+x^{-1}|\vecta)
 =
\begin{cases}
1 &\text{if $d=0$,} \\
2
\big( (x|\vecta)^d - (x^{-1}|\vecta)^d \big)
\dfrac{ x + x^{-1} }
      { x - x^{-1} }
&\text{if $d \ge 1$.}
\end{cases}
$$
Then the factorial symplectic $Q$-function $Q^C_\lambda(\vectx|\vecta)$ 
is obtained from the generalized $P$-function associated with $\GG = \{ g_d(u|\vecta) \}_{d=0}^\infty$ 
by a simple transformation of variables:
$$
Q^C_\lambda(\vectx|\vecta)
 =
P^{\GG}_\lambda(\vectx+\vectx^{-1}).
$$
\end{prop}

In particular, for a single variable $\vectx = (x)$, we have $Q^C_{(r)}(x|\vecta) = \tilde{g}_r(x|\vecta)$.

\begin{corollary}
\label{cor:fQC}
\begin{enumerate}
\item[(1)]
For a strict partition $\lambda$ of length $\le n$, we have
\begin{equation}
\label{eq:Nimmo_fQC}
Q^C_\lambda(\vectx|\vecta)
 =
\begin{cases}
 \dfrac{ 1 }{ \tilde{\Delta}(\vectx) }
 \Pf \begin{pmatrix}
  \tilde{A}(\vectx) & \tilde{W}_\lambda(\vectx|\vecta) \\
  - \trans \tilde{W}_\lambda(\vectx|\vecta) & O
 \end{pmatrix}
 &\text{if $n+l(\lambda)$ is even,}
\\
 \dfrac{ 1 }{ \tilde{\Delta}(\vectx) }
 \Pf \begin{pmatrix}
  \tilde{A}(\vectx) & \tilde{W}_{\lambda^0}(\vectx|\vecta) \\
  - \trans \tilde{W}_{\lambda^0}(\vectx|\vecta) & O
 \end{pmatrix}
 &\text{if $n+l(\lambda)$ is odd,} \\
\end{cases}
\end{equation}
where $\tilde{A}(\vectx)$ and $\tilde{\Delta}(\vectx)$ are given by (\ref{eq:AC}) and 
$\tilde{W}_{(\alpha_1, \dots, \alpha_r)}(\vectx|\vecta)
 = 
\big( \tilde{g}_{\alpha_j}(x_i|\vecta) \big)_{1 \le i \le n, 1 \le j \le r}$.
\item[(2)]
If we adopt (\ref{eq:Nimmo_fQC}) as as the definitions of $Q^C_\lambda(\vectx|\vecta)$ 
for a general strict partition $\lambda$, 
then we have
\begin{equation}
\label{eq:Schur_fQC}
Q^C_\lambda(\vectx|\vecta)
 =
\Pf \Big(
 Q^C_{(\lambda_i,\lambda_j)}(\vectx|\vecta)
\Big)_{1 \le i, j \le m},
\end{equation}
where $m = l(\lambda)$ or $l(\lambda)+1$ according whether $l(\lambda)$ is even or odd.
\end{enumerate}
\end{corollary}

\subsection{
Universal factorial symplectic $Q$-functions
}

In this subsection we lift factorial symplectic $Q$-functions to symmetric functions 
and introduce skew factorial symplectic $Q$-functions.
Let $\Lambda[\vecta] = \Lambda \otimes \Rat[\vecta]$ be the ring of symmetric functions 
in $X = \{ x_1, x_2, \dots \}$ with coefficients in 
$\Rat[\vecta] = \Rat[a_0, a_1, \dots]$, 
and put $\Gamma[\vecta] = \Gamma \otimes \Rat[\vecta] = \Rat[\vecta][q_r : r \ge 1]$.

Since $(x|\vecta)^r = \sum_{k=0}^r e_{r-k}(a_0, \dots, a_{r-1}) x^k$, we have
\begin{equation}
\label{eq:fg_by_g}
\tilde{g}_r(x|\vecta)
 =
\sum_{k=1}^r e_{r-k}(a_0, \dots, a_{r-k}) \tilde{g}_k(x),
\end{equation}
where $e_d$ is the $d$th elementary symmetric polynomial 
and $\tilde{g}_d$ is the Laurent polynomial given by (\ref{eq:g}).
By using this relation, we can express a factorial symplectic $Q$-function
as a linear combination of symplectic $Q$-functions as follows:

\begin{lemma}
\label{lem:fQC_by_QC}
For a strict partition $\lambda$, we have
\begin{equation}
\label{eq:fQC_by_QC}
Q^C_\lambda(\vectx|\vecta)
 =
\sum_\mu d_{\lambda,\mu} Q^C_\mu(\vectx),
\end{equation}
where $\mu$ runs over all strict partitions such that $\mu \subset \lambda$ and $l(\mu) = l(\lambda)$, 
and the coefficients $d_{\lambda,\mu}$ are given by
\begin{equation}
\label{eq:d}
d_{\lambda,\mu}
 =
\det \Big( e_{\lambda_i-\mu_j}(a_0, \dots, a_{\lambda_i-1}) \Big)_{1 \le i, j \le l(\lambda)}.
\end{equation}
\end{lemma}

\begin{demo}{Proof}
We put $l = l(\lambda)$.
First we consider the case where $n+l$ is even.
We start with the Nimmo-type formula (\ref{eq:Nimmo_fQC}) for $Q^C_\lambda(\vectx|\vecta)$ 
and use the multilinearity and alternating property of Pfaffians together with (\ref{eq:fg_by_g}).
Then, by using the Nimmo-type formula (\ref{eq:Nimmo_QC}) for symplectic $Q$-functions, we have
\begin{align*}
Q^C_\lambda(\vectx|\vecta)
 &=
\frac{ 1 }{ \tilde{\Delta}(\vectx) }
\sum_\alpha
 \prod_{i=1}^l e_{\lambda_i - \alpha_i}(a_0, \dots, a_{\lambda_i-1})
\Pf \begin{pmatrix}
 \tilde{A}(\vectx) & \tilde{W}_\alpha(\vectx) \\
 -\trans \tilde{W}_\alpha(\vectx) & O
\end{pmatrix}
\\
 &=
\frac{ 1 }{ \tilde{\Delta}(\vectx) }
\sum_\mu
 \det \Big( e_{\lambda_i - \mu_j}(a_0, \dots, a_{\lambda_i-1}) \Big)_{1 \le i, j \le l}
\Pf \begin{pmatrix}
 \tilde{A}(\vectx) & \tilde{W}_\mu(\vectx) \\
 -\trans \tilde{W}_\mu(\vectx) & O
\end{pmatrix}
\\
 &=
\sum_\mu
 \det \Big( e_{\lambda_i - \mu_j}(a_0, \dots, a_{\lambda_i-1}) \Big)_{1 \le i, j \le l}
Q_\mu(\vectx),
\end{align*}
where $\alpha = (\alpha_1, \dots, \alpha_l)$ runs over all sequence of positive integers, 
$\mu$ runs over all strict partitions such that $l(\mu) = l(\lambda)$.
Also we have $\det \Big( e_{\lambda_i - \mu_j}(a_0, \dots, a_{\lambda_i-1}) \Big)_{1 \le i, j \le l} = 0$ 
unless $\lambda \supset \mu$.
We can prove the case where $n+l$ is odd in a similar way.
\end{demo}

By using this lemma, we can introduce universal factorial symplectic $Q$-functions.

\begin{prop}
\label{prop:ufQC}
For a strict partition $\lambda$, we define the corresponding 
\emph{universal factorial symplectic $Q$-function} 
$\Qfun^C_\lambda(X|\vecta) \in \Gamma[\vecta]$ by
\begin{equation}
\label{eq:ufQC}
\Qfun^C_\lambda(X|\vecta)
 =
\sum_\mu d_{\lambda,\mu} \Qfun^C_\mu(X),
\end{equation}
where $\Qfun^C_\mu$ is the universal symplectic $Q$-function, 
$\mu$ runs over all strict partitions such that $\mu \subset \lambda$ and $l(\mu) = l(\lambda)$, 
and the coefficients $d_{\lambda,\mu}$ are given by (\ref{eq:d}).
Then we have
\begin{enumerate}
\item[(1)]
The symmetric function $\Qfun^C_\lambda(X|\vecta)$ is the unique one satisfying
$$
\tilde{\pi}_n (\Qfun^C_\lambda(X|\vecta))
 = 
Q^C_\lambda(x_1, \dots, x_n|\vecta),
$$
where $\tilde{\pi}_n : \Lambda \to \Rat[x_1^{\pm 1}, \dots, x_n^{\pm 1}]^{W_n}$ is the ring homomorphism 
given by (\ref{eq:pi}).
\item[(2)]
For a strict partition $\lambda$, we have
$$
\Qfun^C_\lambda(X|\vecta)
 =
\Pf \Big(
 \Qfun^C_{(\lambda_i,\lambda_j)}(\vectx|\vecta)
\Big)_{1 \le i, j \le m},
$$
where $m = l(\lambda)$ or $l(\lambda)+1$ according whether $l(\lambda)$ is even or odd.
\item[(3)]
The universal factorial symplectic $Q$-functions $\{ \Qfun^C_\lambda(X|\vecta) : \lambda \in \SPar \}$ 
form a $\Rat[\vecta]$-basis of $\Gamma[\vecta]$.
\end{enumerate}
\end{prop}

\begin{demo}{Proof}
By using Lemma~\ref{lem:proj}, 
(1) and (2) follows from Lemma~\ref{lem:fQC_by_QC} and (\ref{eq:Schur_fQC}) respectively.
Since $\Qfun^C_\lambda(X|\vecta) = \Qfun^C_\lambda(X) + \text{lower degree terms}$ by definition, 
we can derive the claim (3) from Proposition~\ref{prop:uQC_basis} (2).
\end{demo}

\begin{remark}
By \cite[Theorem~3.5${}^*$]{Wachs85}, the coefficient $d_{\lambda,\mu}$ given by (\ref{eq:d}) 
is equal to the multivariate generating function of 
row-strict tableaux of shape $\lambda/\mu$ such that entries of the $i$th row are bounded by $\lambda_i-1$.
Here a row-strict tableau of shape $\lambda/\mu$ is 
a filling of the boxes of the skew shifted diagram $S(\lambda/\mu)$ 
with nonnegative integers such that each row is strictly increasing and each column is weakly increasing.
\end{remark}

Now we introduce skew factorial symplectic $Q$-functions.
By the definition (\ref{eq:fQC_by_QC}), we see that $\Qfun^C_{(r)}(X|\vecta)$ depends on $a_0, \dots, a_{r-1}$, 
so we write $\Qfun^C_{(r)}(X|a_0, \dots, a_{r-1})$.
For two nonnegative integer $r$ and $k$, 
we define $\Rfun^C_{r/k}(X|\vecta) \in \Lambda[\vecta]$ by
\begin{equation}
\label{eq:R}
\Rfun^C_{r/k}(X|\vecta)
 =
\begin{cases}
 \Qfun^C_{(r)}(X|a_0, a_1, \dots, a_{r-1})
 &\text{if $k=0$,} \\
 \Qfun^C_{(r-k)}(X|0, a_{k+1}, \dots, a_{r-1})
 &\text{if $1 \le k \le r-1$,} \\
 1 &\text{if $k=r$,} \\
 0 &\text{otherwise.}
\end{cases}
\end{equation}
For sequences of nonnegative integers $\alpha = (\alpha_1, \dots, \alpha_r)$ 
and $\beta = (\beta_1, \dots, \beta_s)$, we put
$$
\tilde{K}_\alpha(X|\vecta)
 =
\Big( \Qfun^C_{(\alpha_i,\alpha_j)}(X|\vecta) \Big)_{1 \le i, j \le r},
\quad
\tilde{M}_{\alpha/\beta}(X|\vecta)
 =
\Big( \Rfun^C_{\alpha_i/\beta_{s+1-j}}(X|\vecta) \Big)_{1 \le i \le r, 1 \le j \le s}.
$$
Given two strict partition $\lambda$ and $\mu$, we define 
\begin{equation}
\label{eq:fskewQC}
\Qfun^C_{\lambda/\mu}(X|\vecta)
 =
\begin{cases}
 \Pf \begin{pmatrix}
  \tilde{K}_\lambda(X|\vecta) & \tilde{M}_{\lambda/\mu}(X|\vecta) \\
  -\trans \tilde{M}_{\lambda/\mu}(X|\vecta) & O
 \end{pmatrix}
&\text{if $l(\lambda)+l(\mu)$ is even,}
\\
 \Pf \begin{pmatrix}
  \tilde{K}_\lambda(X|\vecta) & \tilde{M}_{\lambda/\mu^0}(X|\vecta) \\
  -\trans \tilde{M}_{\lambda/\mu^0}(X|\vecta) & O
 \end{pmatrix}
&\text{if $l(\lambda)+l(\mu)$ is odd,}
\end{cases}
\end{equation}
where $\mu^0 = (\mu_1, \dots, \mu_{l(\mu)},0)$.
We call $\Qfun^C_{\lambda/\mu}(X|\vecta)$ the \emph{universal skew factorial symplectic $Q$-function}.
For a finite number of variables $\vectx = (x_1, \dots, x_n)$, we put
$$
R_{r/k}(\vectx|\vecta)
 =
\tilde{\pi}_n (\Rfun_{r/k}(X|\vecta)),
\quad
Q^C_{\lambda/\mu}(\vectx|\vecta)
 =
\tilde{\pi}_n (\Qfun^C_{\lambda/\mu}(X|\vecta)).
$$

By comparing (\ref{eq:Schur_fQC}) with(\ref{eq:fskewQC}), 
we see that $\Qfun^C_{\lambda/\emptyset}(X|\vecta) = \Qfun^C_\lambda(X|\vecta)$.
Also by the same argument as in the proof of \cite[Proposition~4.4]{Okada20}, we have

\begin{prop}
\label{prop:ufskewQC}
For two strict partitions $\lambda$ and $\mu$, 
we have $\Qfun^C_{\lambda/\mu}(X|\vecta) = 0$ unless $\lambda \supset \mu$.
\end{prop}

The following is the key property of skew symplectic $Q$-functions, 
and play an important role in the proof of Theorems~\ref{thm:fac_tableaux} 
and \ref{thm:fac_skew_tableaux}.

\begin{prop}
\label{prop:fac_sep_var}
Let $X$ and $Y$ be two disjoint sets of infinitely many variables.
For strict partitions $\lambda$ and $\nu$, we have
\begin{gather}
\label{eq:fac_sep_var1}
\Qfun^C_\lambda(X \cup Y|\vecta)
 =
\sum_\mu
\Qfun^C_{\lambda/\mu}(X|\vecta)
\Qfun^C_\mu(Y|\vecta),
\\
\label{eq:fac_sep_var2}
\Qfun^C_{\lambda/\nu}(X \cup Y|\vecta)
 =
\sum_\mu
\Qfun^C_{\lambda/\mu}(X|\vecta)
\Qfun^C_{\mu/\nu}(Y|\vecta).
\end{gather}
\end{prop}

In order to prove this proposition, we need the following lemma.

\begin{lemma}
\label{lem:Q=Rg}
For variables $x_1, \dots, x_n$ and $y$, we have
\begin{equation}
\label{eq:Q=Rg}
Q^C_{(r)}(x_1, \dots, x_n, y |\vecta)
 =
\sum_{k=0}^r
 R^C_{r/k}(x_1, \dots, x_n|\vecta) \tilde{g}_k(y|\vecta).
\end{equation}
\end{lemma}

\begin{demo}{Proof}
We express the both sides of (\ref{eq:Q=Rg}) as linear combinations of $Q^C_{(i)}(\vectx) \tilde{g}_j(y)$.
By using (\ref{eq:fQC_by_QC}), (\ref{eq:sep_var_length1}) and (\ref{eq:n=1}), we have
$$
Q^C_{(r)}(\vectx,y|\vecta)
 =
\sum_{s=1}^r e_{r-k}(a_0, \dots, a_{r-1}) Q^C_{(s)}(\vectx,y)
 =
\sum_{s=1}^r e_{r-s}(a_0, \dots, a_{r-1}) \sum_{j=0}^s Q^C_{(s-j)}(\vectx) \tilde{g}_j(y).
$$
Similarly we have
\begin{align*}
&
\sum_{k=0}^r
 R^C_{r/k}(x_1, \dots, x_n|\vecta) \tilde{g}_k(y)
\\
&\quad
=
\sum_{i=1}^r e_{r-i}(a_0, \dots, a_{r-1}) Q^C_{(i)}(\vectx)
\\
&\quad\quad
+
\sum_{k=1}^{r-1} 
 \sum_{i=1}^{r-k} e_{r-k-i}(0, a_{k+1}, \dots, a_{r-1}) Q^C_{(i)}(\vectx)
 \sum_{j=1}^k e_{k-j}(a_0, \dots, a_{k-1}) \tilde{g}_j(y)
\\
&\quad\quad
+
\sum_{j=1}^r e_{r-j}(a_0, \dots, a_{r-1}) \tilde{g}_j(y).
\end{align*}
Now, by using the relation (\cite[Lemma~6.4 (1)]{Okada20})
\begin{equation}
\label{eq:rel-e}
\sum_{k=j}^{r-i}
 e_{k-j}(a_0, \dots, a_{k-1}) e_{r-k-i}(a_{k+1}, \dots, a_{r-1})
 =
e_{r-i-j}(a_0, \dots, a_{r-1}),
\end{equation}
we obtain the desired identity (\ref{eq:Q=Rg}).
\end{demo}

\begin{demo}{Proof of Proposition~\ref{prop:fac_sep_var}}
Once the relation (\ref{eq:Q=Rg}) is established, we can use the same argument as in the proof of 
\cite[Theorem~4.2 and Proposition~4.5]{Okada20} to show
$$
Q^C_\lambda(x_1, \dots, x_n, y_1, \dots, y_m)
 =
\sum_\mu Q_{\lambda/\mu}(x_1, \dots, x_n) Q_\mu(y_1, \dots, y_m),
$$
where $\mu$ runs over all strict partitions.
Then by using Lemma~\ref{lem:proj} (2), we obtain (\ref{eq:fac_sep_var1}).
We can derive (\ref{eq:fac_sep_var2}) from (\ref{eq:fac_sep_var1}) 
by a standard argument used in the proof of Proposition~\ref{prop:sep_var}.
\end{demo}

\subsection{%
Tableau description for skew factorial symplectic $Q$-funtions
}

Now we can state and prove a skew version of Theorem~\ref{thm:fac_tableaux}.

\begin{theorem}
\label{thm:fac_skew_tableaux}
Let $\vectx = (x_1, \dots, x_n)$ be a sequence of $n$ indeterminates.
Assume $a_0 = 0$.
For strict partitions $\lambda$ and $\mu$, we have
\begin{equation}
\label{eq:fac_skew_tableaux}
Q^C_{\lambda/\mu}(\vectx|\vecta)
 =
\sum_{T \in \QTab^C_n(\lambda/\mu)} (\vectx|\vecta)^T.
\end{equation}
\end{theorem}

The method of the proof of Theorem~\ref{thm:fac_skew_tableaux} is exactly the same as 
the proof of Theorem~\ref{thm:tableaux} for skew symplectic $Q$-functions.
The following lemma can be proved in a similar way to the proof of Lemma~\ref{lem:algQ} 
by using properties of factorial symplectic $Q$-functions.

\begin{lemma}
\label{lem:fac_algQ}
Let $\lambda$ and $\mu$ be strict partitions such that $\lambda \supset \mu$.
\begin{enumerate}
\item[(1)]
We have
$$
Q^C_{\lambda/\mu}(x_1, \dots, x_n|\vecta)
 =
\sum \prod_{i=1}^n Q^C_{\mu^{(i)}/\mu^{(i-1)}} (x_i|\vecta),
$$
where the sum is taken over all sequences $\mu = \mu^{(0)} \subset \mu^{(1)} \subset \dots \subset 
\mu^{(n-1)} \subset \mu^{(n)} = \lambda$.
\item[(2)]
For a single variable $x_1$, we have $Q^C_{\lambda/\mu}(x_1|\vecta) = 0$ unless $l(\lambda) - l(\mu) \le 1$.
\item[(3)]
If $l(\lambda) - l(\mu) \le 1$, then we have
$$
Q^C_{\lambda/\mu}(x_1|\vecta)
 =
\det \Big( R^C_{\lambda_i/\mu_j}(x_1|\vecta) \Big)_{1 \le i, j \le l(\lambda)}.
$$
\end{enumerate}
\end{lemma}

On the combinatorial side, we put
$$
Q^{\text{tab}}_{\lambda/\mu}(\vectx|\vecta)
 =
\sum_{T \in \QTab^C_n(\lambda/\mu)} (\vectx|\vecta)^T
$$
and prove they satisfy the same relations as $Q^C_{\lambda/\mu}(\vectx|\vecta)$.

\begin{lemma}
\label{lem:fac_combQ}
Let $\lambda$ and $\mu$ be strict partitions such that $\lambda \supset \mu$.
\begin{enumerate}
\item[(1)]
We have
$$
Q^{\text{tab}}_{\lambda/\mu}(x_1, \dots, x_n|\vecta)
 =
\sum \prod_{i=1}^n Q^{\text{tab}}_{\mu^{(i)}/\mu^{(i-1)}} (x_i|\vecta),
$$
where the sum is taken over all sequences $\mu = \mu^{(0)} \subset \mu^{(1)} \subset \dots \subset 
\mu^{(n-1)} \subset \mu^{(n)} = \lambda$.
\item[(2)]
For a single variable $x_1$, 
we have $Q^{\text{tab}}_{\lambda/\mu}(x_1|\vecta) = 0$ unless $l(\lambda) - l(\mu) \le 1$.
\item[(3)]
If $l(\lambda) - l(\mu) \le 1$, then we have
$$
Q^{\text{tab}}_{\lambda/\mu}(x_1|\vecta)
 =
\det \Big( Q^{\text{tab}}_{(\lambda_i)/(\mu_j)}(x_1|\vecta) \Big)_{1 \le i, j \le l(\lambda)},
$$
where $Q^{\text{tab}}_{(r)/(k)}(x) = 0$ for $r<k$.
\end{enumerate}
\end{lemma}

\begin{demo}{Proof}
We give a proof of (3).
We consider the same graph $G$ (see Figure~\ref{fig:G}) 
used in the proof of Theorem~\ref{thm:tableaux} with different edge weight given by
\begin{gather*}
\wt(A_i,B_i) = \wt (B_i,C_i) = 1,
\\
\wt(A_i,B_{i+1}) = x - a_i,
\quad
\wt(B_i,B_{i+1}) = x + a_i,
\\
\wt(B_i, C_{i+1}) = x^{-1} - a_i,
\quad
\wt(C_i, C_{i+1}) = x^{-1} + a_i.
\end{gather*}
Then by applying the Lindstr\"om--Gessel--Viennot lemma, we obtain
$$
Q^{\text{tab}}_{\lambda/\mu}(x_1|\vecta)
 =
\det \Big( Q^{\text{tab}}_{(\lambda_i)/(\mu_j)}(x_1|\vecta) \Big)_{1 \le i, j \le l(\lambda)}.
$$
\end{demo}

Now we can finish the proof of Theorem~\ref{thm:fac_skew_tableaux}.

\begin{demo}{Proof of Theorem~\ref{thm:fac_skew_tableaux}}
By comparing Lemma~\ref{lem:fac_algQ} with Lemma~\ref{lem:fac_combQ}, 
it is enough to show that $R^C_{r/k}(x|\vecta) = Q^{\text{tab}}_{(r)/(k)}(x|\vecta)$.
We put $s = r-k$ and compute $Q^{\text{tab}}_{(r)/(k)}(x|\vecta)$ explicitly 
by listing all possible primed shifted tableaux of shape $(r)/(k)$.
Then we have
\begin{align*}
&
Q^{\text{tab}}_{(r)/(k)}(x|\vecta)
\\
&\quad
=
(x-a_k) \prod_{i=1}^{s-1} (x+a_{k+i})
+
(x+a_k) \prod_{i=1}^{s-1} (x+a_{k+i})
\\
&\quad\quad
+
\sum_{l=1}^{s-1}
 (x-a_k) \prod_{i=1}^{l-1} (x + a_{k+i}) \cdot (x^{-1}-a_{k+l}) \prod_{i=l+l}^{s-1} (x^{-1} + a_{k+i})
\\
&\quad\quad
+
\sum_{l=1}^{s-1}
 (x+a_k) \prod_{i=1}^{l-1} (x + a_{k+i}) \cdot (x^{-1}-a_{k+l}) \prod_{i=l+l}^{s-1} (x^{-1} + a_{k+i})
\\
&\quad\quad
+
\sum_{l=1}^{s-1}
 (x-a_k) \prod_{i=1}^{l-1} (x + a_{k+i}) \cdot (x^{-1}+a_{k+l}) \prod_{i=l+l}^{s-1} (x^{-1} + a_{k+i})
\\
&\quad\quad
+
\sum_{l=1}^{s-1}
 (x+a_k) \prod_{i=1}^{l-1} (x + a_{k+i}) \cdot (x^{-1}+a_{k+l}) \prod_{i=l+l}^{s-1} (x^{-1} + a_{k+i})
\\
&\quad\quad
+
(x^{-1}-a_k) \prod_{i=1}^{s-1} (x^{-1} + a_{k+i})
+
(x^{-1}+a_k) \prod_{i=1}^{s-1} (x^{-1} + a_{k+i})
\\
&\quad
=
2 x \prod_{i=1}^{s-1} (x+a_{k+i})
 +
4 \sum_{l=1}^{s-1}
 \prod_{i=1}^{l-1} (x + a_{k+i})
 \prod_{i=l+l}^{s-1} (x^{-1} + a_{k+i})
 +
2 x^{-1} \prod_{i=1}^{s-1} (x^{-1} + a_{k+i})
\\
&\quad
=
2 x \sum_{l=0}^{s-1} e_{s-1-l}(a_{k+1}, \dots, a_{r-1}) x^l
\\
&\quad\quad
+
4 \sum_{l=1}^{s-1}
 \sum_{i=0}^{l-1} e_{l-1-i}(a_{k+1}, \dots, a_{k+l-1}) x^i
 \sum_{j=0}^{s-l-1} e_{s-l-1-j}(a_{k+l+1}, \dots, a_{r-1}) x^{-j}
\\
&\quad\quad
+
2 x^{-1} \sum_{l=0}^{s-1} e_{s-1-l}(a_{k+1}, \dots, a_{r-1}) x^{-l}.
\end{align*}
Here we appeal to the relation (\ref{eq:rel-e}), which is used in the proof of Lemma~\ref{lem:Q=Rg}.
Then by using (\ref{eq:fg_by_g}) and (\ref{eq:R}), we obtain
\begin{align*}
Q^{\text{tab}}_{(r)/(k)}(x|\vecta)
&=
2 \sum_{l=1}^s
 e_{s-l}(a_{k+1}, \dots, a_{r-1})
 \big( x^l + 2 x^{l-2} + 2 x^{l-4} + \dots + 2 x^{-l+2} + x^{-l} \big)
\\
&=
\sum_{l=1}^s
 e_{s-l}(0, a_{k+1}, \dots, a_{r-1}) \tilde{g}_l(x)
 =
\tilde{g}_l(x|0, a_{k+1}, \dots, a_{r-1})
 =
R_{r/k}(x|\vecta).
\end{align*}
This completes the proof of Theorem~\ref{thm:fac_skew_tableaux}.
\end{demo}



\begin{thebibliography}{WW}

\bibitem{Cho13}
S.~Cho,
A new Littlewood--Richardson rule for Schur $P$-functions,
Trans. Amer. Math. Soc. {\bfseries 365} (2013), 939--972.

\bibitem{FK18}
A.~M.~Foley and R.~C.~King,
Factorial $Q$-functions and Tokuyama identities for classical Lie groups, 
European J. Combin. {\bfseries 73} (2018), 89--113.

\bibitem{GJKKK14}
D.~Grantcharov, J.~H.~Jung, S.-J.~Kang, M.~Kashiwara, and M.~H.~Kim,
Crystal bases for the quantum queer superalgebra and semistandard decomposition tableaux,
Trans. Amer. Math. Soc. {\bfseries 366} (2014), 457--489.

\bibitem{HK07}
A.~M.~Hamel and R.~C.~King,
Bijective proofs of shifted tableau and alternating sign matrix identities, 
J. Algebr. Comb. {\bfseries 25} (2007), 417--458.

\bibitem{King76}
R.~C.~King,
Weight multiplicities for the classical groups,
in ``Group Theoretical Methods in Physics (Fourth Internat. Colloq., Nijmegen, 1975)'', 
Lecture Notes in Phys. {\bfseries 50}, Springer, Berlin, 1976, 
pp. 490--499. 

\bibitem{King90}
R.~C.~King,
$S$-functions and characters of Lie algebras and superalgebras,
in ``Invariant theory and tableaux (Minneapolis, MN, 1988)'',
IMA Vol. Math. Appl. {\bfseries 19}, Springer, New York, 1990,
pp. 226--261.

\bibitem{KH07}
R.~C.~King and A.~M.~Hamel,
Combinatorial realisation of Hall--Littlewood polynomials at $t = -1$,
Proceedings of the 19th International Conference on Formal Power Series and Algebraic Combinatorics 
(Tianjin, July 2--6, 2007),
available at {\tt http://igm.univ-mlv.fr/{\textasciitilde}fpsac/FPSAC07/SITE07/PDF-Proceedings/Posters/75.pdf}

\bibitem{Koike89}
K.~Koike,
On the decomposition of tensor products of the representations of the classical groups:
By means of the universal characters,
Adv. Math. {\bfseries 74} (1989), 57--86.

\bibitem{KT87}
K.~Koike and I.~Terada,
Young-diagrammatic methods for the representation theory of the classical groups of type $B_n$, $C_n$, $D_n$, 
J. Algebra {\bfseries 107} (1987), 466--511.

\bibitem{Macdonald95}
I.~G.~Macdonald,
``Symmetric Functions and Hall Polynomials, 2nd edition'',
Oxford Univ. Press, 1995.

\bibitem{Macdonald00}
I.~G.~Macdonald,
Orthogonal polynomials associated with root systems,
S\'em. Lothar. Combin. {\bfseries 45} (2000/01), Art. B45a, 40pp.

\bibitem{Morris64}
A.~O.~Morris,
A note on the multiplication of Hall functions,
J. London Math. Soc. {\bfseries 39} (1964), 481--488.

\bibitem{Nimmo90}
J.~J.~C.~Nimmo,
Hall--Littlewood symmetric functions and the BKP equation,
J. Phys. A {\bfseries 23} (1990), 751--760.

\bibitem{Okada14}
S.~Okada,
Schur-type Pfaffians and their applications to symmetric function identities,
in ``Abstract of the Special Session on Infinite Analysis, MSJ 2014 Autumn Meeting'', 
Mathematical Society of Japan, 2014, 
pp.21--31 (in Japanese).

\bibitem{Okada19}
S.~Okada,
Pfaffian formulas and Schur $Q$-function identities,
Adv. Math. {\bfseries 353} (2019), 446--470, 

\bibitem{Okada20}
S.~Okada, 
A generalization of Schur's $P$- and $Q$-functions,
S\'em. Lothar. Combin. {\bfseries 81} (2019/20), Art. B81k, 50pp.

\bibitem{PJ91}
P.~Pragacz and T.~J\'ozefiak,
A determinantal formula for skew $Q$-functions, 
J. London Math. Soc. (2) {\bf 43} (1991), 76--90.

\bibitem{Sagan87}
B.~E.~Sagan,
Shifted tableaux, Schur $Q$-functions, and a conjecture of R.~Stanley,
J. Combin. Theory Ser. A {\bfseries 45} (1987), 62--103.

\bibitem{Schur}
I.~Schur,
\"Uber die Darstellung der symmetrischen und der alternierenden Gruppe durch gebrochene lineare Substitutionen, 
J. Reine Angew. Math. {\bfseries 139} (1911), 155--250.

\bibitem{Stembridge89}
J.~R.~Stembridge,
Shifted tableaux and the projective representations of symmetric groups,
Adv. Math. {\bfseries 74} (1989), 87--134.

\bibitem{SF}
J.~R.~Stembridge,
The {\tt SF} package for symmetric functions, 
available at {\tt http://www.math.lsa.umich.edu/{\textasciitilde}jrs/maple.html\#SF}

\bibitem{Sundaram90}
S.~Sundaram,
The Cauchy identity for $Sp(2n)$,
J. Combin. Theory Ser. A {\bf 53} (1990), 209--238.

\bibitem{Wachs85}
M.~L.~Wachs,
Flagged Schur functions, Schubert polynomials, and symmetrizing operators, 
J. Combin. Theory Ser. A {\bfseries 40} (1985), 276--289.

\bibitem{Worley84}
D.~R.~Worley,
``A Theory of Shifted Young Tableaux'',
Ph.D. Thesis, 
Massachusetts Institute of Technology, 1984.

\end{thebibliography}
\end{document}